\documentclass[11pt]{amsart}
\usepackage{amsmath}
\usepackage{amsxtra}
\usepackage{amscd}
\usepackage{amsthm}
\usepackage{amsfonts}
\usepackage{amssymb}
\usepackage{mathrsfs}
\usepackage[all]{xy}
\usepackage[pdftex]{color,graphicx}

\newtheorem{prop}{Proposition}[section]
\newtheorem{thm}[prop]{Theorem}
\newtheorem{lem}[prop]{Lemma}
\newtheorem{cor}[prop]{Corollary}

\theoremstyle{definition}

\newtheorem{defn}[prop]{Definition}
\newtheorem{rem}[prop]{Remark}
\newtheorem{exam}[prop]{Example}

\newcommand{\lbl}[1]{ \label{#1}}
\newcommand{\editorial}[1]{}

\newcommand{\uG}{\underline{G}}
\newcommand{\uH}{\underline{H}}

\newcommand{\al}{\alpha}
\newcommand{\be}{\beta}
\newcommand{\ga}{\gamma}

\newcommand{\ep}{\epsilon}
\newcommand{\ze}{\zeta}

\newcommand{\la}{\lambda}
\newcommand{\om}{\omega}

\newcommand{\Ga}{\Gamma}
\newcommand{\De}{\Delta}

\newcommand{\Om}{\Omega}

\newcommand{\cA}{\mathcal A}
\newcommand{\cB}{\mathcal B}

\newcommand{\cD}{\mathcal D}
\newcommand{\cE}{\mathcal E}
\newcommand{\cF}{\mathcal F}
\newcommand{\cG}{\mathcal G}

\newcommand{\cL}{\mathcal L}

\newcommand{\cN}{\mathcal N}
\newcommand{\cO}{\mathcal O}
\newcommand{\cP}{\mathcal P}

\newcommand{\cT}{\mathcal T}

\newcommand{\cV}{\mathcal V}
\newcommand{\cW}{\mathcal W}

\newcommand{\bA}{\mathbb A}

\newcommand{\bC}{\mathbb C}

\newcommand{\bF}{\mathbb F}
\newcommand{\bG}{\mathbb G}

\newcommand{\bK}{\mathbb K}

\newcommand{\bM}{\mathbb M}
\newcommand{\bN}{\mathbb N}
\newcommand{\bO}{\mathbb O}
\newcommand{\bP}{\mathbb P}
\newcommand{\bQ}{\mathbb Q}
\newcommand{\bR}{\mathbb R}

\newcommand{\bZ}{\mathbb Z}

\newcommand{\fg}{\mathfrak g}

\newcommand{\fl}{\mathfrak l}
\newcommand{\fp}{\mathfrak p}
\newcommand{\fq}{\mathfrak q}
\newcommand{\fr}{\mathfrak r}

\newcommand{\sK}{\mathscr K}
\newcommand{\sL}{\mathscr L}
\newcommand{\sN}{\mathscr N}
\newcommand{\sU}{\mathscr U}

\newcommand{\qq}{\quad}
\newcommand{\QQ}{\qquad\quad}

\newcommand{\bmat}{\left ( \begin{matrix} }
\newcommand{\emat}{\end{matrix} \right ) }

\DeclareMathOperator{\val}{val}
\DeclareMathOperator{\Irr}{Irr}
\DeclareMathOperator{\Spec}{Spec}
\DeclareMathOperator{\Ind}{Ind}
\DeclareMathOperator{\Res}{Res}
\DeclareMathOperator{\Stab}{Stab}
\DeclareMathOperator{\GL}{GL}
\DeclareMathOperator{\M}{M}
\DeclareMathOperator{\Id}{Id}
\DeclareMathOperator{\Grass}{Grass}
\DeclareMathOperator{\Ad}{Ad}
\DeclareMathOperator{\Lie}{Lie}
\DeclareMathOperator{\Aut}{Aut}
\DeclareMathOperator{\trace}{trace}

\newcommand{\X}{\mathscr{X}} 
\newcommand{\Y}{\mathscr{Y}} 
\DeclareMathOperator{\ac}{ac} 
\newcommand{\eff}{\mathfrak{f}}
\title{Arithmetic Groups Have Rational Representation Growth}

\author{Nir Avni}

\begin{document}

\maketitle

\begin{abstract} Let $\Ga$ be an arithmetic lattice in a semisimple algebraic group over a number field. We show that if $\Ga$ has the congruence subgroup property, then the number of $n$-dimensional irreducible representations of $\Ga$ grows like $n^\al$, where $\al$ is a rational number.
\end{abstract}

\section{Introduction}

\subsection{Representation Zeta Functions} This article is concerned with counting the number of representations of arithmetic groups. Suppose that $\Ga$ is a finitely generated group, and assume that $\Ga$ has finitely many irreducible complex representation of any fixed dimension, up to equivalence. Denote the number of irreducible complex representations of $\Ga$ of dimension $n$, up to equivalence by $r_n(\Ga)$. In \cite{LM}, the sequence $r_n(\Ga )$ is called the {\em representation growth sequence of $\Ga$}. If the sequence $r_n(\Ga)$ is bounded by a polynomial in $n$, then it is useful to consider the following generating function:

\begin{defn} \lbl{defn:rep.zeta.func} The representation zeta function of $\Ga$ is the following function of $s\in\bC$:
\[
\ze _{\Ga }(s)=\sum _{n=1} ^{\infty} r_n (\Ga )n^{-s}  = \sum _{\rho \in \Irr \Ga} (\dim\rho )^{-s},
\]
where $\Irr\Ga$ denotes the set of finite dimensional, complex, and irreducible representations of $\Ga$.
\end{defn}

Note that if the sequence $r_n (\Ga )$ grows polynomially, then the series above converges in some half plane of the form $\{ s | \Re(s)>\al \}$. The infimum of the set of $\al \in \bR$ such that the series in Definition \ref{defn:rep.zeta.func} converges absolutely at $s=\al$, is called the {\em abscissa of convergence of $\ze _{\Ga}(s)$} (or of $\Ga$); We will denote it by $\al _{\Ga }$. The abscissa of converges is related to the rate of growth of the sequence $r_n(\Ga )$ by

\[
\al _{\Ga} = \limsup \limits _{N\to \infty} \frac{\log (r_1 (\Ga )+ \dots + r_N(\Ga ))}{\log N} .
\]

\subsection{Arithmetic Lattices} The groups which we consider in this paper are arithmetic lattices in semisimple algebraic groups over fields of characteristics 0. We remind the reader the construction of such groups. Let $\bK$ be a finite extension of the field of rational numbers $\bQ$. Denote the ring of integers of $\bK$ by $\bO$. For a valuation $v$ of $\bK$, we denote the completion of $\bK$ with respect to the valuation $v$ by $\bK _v$, and we denote the valuation ring of $\bK _v$ by $\bO _v$. Suppose that $\Sigma$ is a finite set of valuations of $\bK$, containing all infinite (i.e. archimedian) valuations.  The ring of $\Sigma$-integers of $\bK$ is the set 
\[
\bO _\Sigma =\{ x\in \bK | (\forall v\not \in \Sigma) \quad v(x)\geq 0\}.
\]

Let $\uG\subset\underline{GL_N}_{\bO_\Sigma}$ be a linear algebraic group scheme over $\Spec \bO _\Sigma$\footnote{Put more simply, we are given a set of polynomials $f_1 , \ldots ,f_k$ in $N^2$ variables, such that the coefficients of the $f_j$'s are in $\bO _\Sigma$, and such that for every ring $R$ and homomorphism $\varphi :\bO _\Sigma \to R$, the set of solutions of the system of equations $(\varphi f_1)(x)=\ldots =(\varphi f_k)(x)=0$ in $R^{N^2}=M_N(R)$ is a subgroup of $GL_N(R)$. We call this set of solutions the {\em $R$ points} of $\uG$ and denote it by $\uG (R)$.} whose generic fiber\footnote{The generic fiber of $\uG$ is the algebraic group $\uG (\overline{\bK})$, where $\overline{\bK}$ is the algebraic closure of $\bK$.} is semisimple, simply connected, and connected. Assume, moreover, that for every non-archimedian valuation $v\in \Sigma$, the group $\uG (\bK _v )$ is non compact. The group $\Ga =\uG (\bO _\Sigma)$ is the arithmetic lattice. It is indeed a lattice, i.e. a discrete subgroup of finite covolume, in the topological group $\prod _{v\in \Sigma} \uG (\bK _{v})$.

Denote the profinite completion of $\Ga$ by $\widehat{\Ga}$, and, similarly, let $\widehat{\bO _\Sigma}$ be the profinite completion of the ring $\bO _\Sigma$. By the Chinese remainder theorem, $\widehat{\bO _\Sigma}=\prod _{v\not \in \Sigma} \bO _v$. We say that $\Ga$ has the {\em congruence subgroup property}, if the kernel of the natural map
\[
\widehat{\Ga }=\widehat{\uG (\bO _\Sigma)} \longrightarrow \uG (\widehat{\bO _\Sigma}) = \prod _{v\not \in \Sigma} \uG (\bO _v)
\]
is finite.

It is known that ``most'' lattices in Lie groups of rank $\geq 2$ have the congruence subgroup property, and a conjecture of Serre asserts that all of them do. See \cite{Ra} for a survey on the congruence subgroup property.

\subsection{Main Theorem} 
In \cite{LM} it was proved that an arithmetic lattice in characteristics 0 has the congruence subgroup property if and only if the sequence $r_n(\Ga )$ grows polynomially. Equivalently, such a lattice, $\Ga$, has the congruence subgroup property if and only if the abscissa of convergence of $\Ga$ is finite. The main result in this paper is the following:

\begin{thm} \lbl{thm:rat.absc} Let $\Ga$ be an arithmetic lattice in characteristics 0 that satisfies the congruence subgroup property. Then $\al_\Ga$---the abscissa of convergence of $\ze _{\Ga}(s)$---is a rational number.
\end{thm}

\begin{rem} If $\Ga$ does not satisfy the congruence subgroup property, then the sequence $r_n(\Ga)$ grows super-polynomially by \cite{LM}, and so the abscissa of convergence of $\ze_\Ga(s)$ is $\infty$.
\end{rem}

Unfortunately, the proof of this theorem does not give a hint about the actual value of the abscissa of convergence of $\Ga$, and, in fact, this value is known only in some very special cases, see \cite{LL} and \cite{AO}.

In the rest of this subsection we describe the method of proof of Theorem \ref{thm:rat.absc}. The proof follows a general strategy of Igusa and Denef, see also \cite{dSG}. If $\Ga=\uG(\bO_\Sigma)$ is an arithmetic lattice that has the congruence subgroup property, then there is a finite index subgroup $\De$ of $\Ga$, such that the representation zeta function of $\De$ has an Euler-like factorization:
\[
\ze_\De(s)=\ze_\infty(s)\times\prod_\fp\ze_\fp(s),
\]
where the product is over all primes of the ring $\bO_\Sigma$, and the local zeta functions $\ze_{\infty}(s)$ and $\ze_\fp(s)$ will be described in Section \ref{sec:Euler.Factorization}. This fact was established in \cite{LL}, and is a consequence of Margulis' super-rigidity theorem. We shall show that the abscissa of convergence is unchanged when passing to a finite-index subgroup. Hence, it is enough to show that the abscissa of convergence of $\De$ is rational. The archimedian local zeta function $\ze_\infty(s)$ was studied in \cite{LL}, where it was shown that it has a rational abscissa of convergence. In order to show that the infinite product $\prod_\fp\ze_\fp(s)$ has rational abscissa of convergence, we will study the dependence of $\ze_\fp(s)$ on the prime ideal $\fp$.

Let $q=|\bO_\Sigma/\fp|$. In contrast to the case considered in \cite{dSG}, the local zeta functions are not rational functions in $q^{-s}$, but rather are of the form
\begin{equation} \lbl{eq:local.zeta.function}
\sum_{i=1}^{N(\fp)}n_i(\fp)^{-s}\cdot f_i(\fp,q^{-s})
\end{equation}
where $f_i(\fp,x)$ are rational functions in $x$ (this is proved in \cite{Jai}).

A sequence of numbers, $k(\fp)$, indexed by the primes of $\bO_\Sigma$, is called geometric, if there is a variety $\sK$, defined over $\bO_\Sigma$, such that for every $\fp$, $k(\fp)$ is equal to the number of points of the variety $\sK$ over the finite field $\bO_\Sigma/\fp$ (this terminology is taken from \cite{Ka}). A reasonable guess is that the numbers $N(\fp),n_i(\fp)$, and the coefficients of the rational functions $f_i(\fp,x)$ that appear in (\ref{eq:local.zeta.function}) are geometric. We make two changes in order to prove this. The first is that we allow $\sK$ to be a definable set, rather than a variety; the second is that we replace $\ze_\fp(s)$ by another sequence, $\xi_\fp(s)$, such that the abscissae of convergence of $\prod_\fp\ze_\fp(s)$ and $\prod_\fp\xi_\fp(s)$ are equal, and then show that $\xi_\fp(s)$ has a geometric formula.

After showing the geometric nature of the ``new'' local zeta functions, we use standard results in Algebraic Number Theory (the Lang-Weil estimates and Chebotarev Density Theorem) to finish the proof of Theorem \ref{thm:rat.absc}.

\subsection{Organization of the Paper}
In Section \ref{sec:Euler.Factorization} we set some notations and review the Euler factorization of representation zeta functions of arithmetic lattices. Section \ref{sec:Algebraic.Preliminaries} is a collection of facts we need from Representation Theory, Algebraic Geometry, and the theory of finite groups. In Section \ref{sec:Definable.Families} we collect necessary facts from the model theory of fields, pseudo-finite fields, and valued fields. In Section \ref{sec:Definable.Families} we also define the notion of $V$-function, which is our main technical tool. In Section \ref{sec:Uniformity.I} we show that local zeta functions (or, rather, approximations of which) are integrals of the same $V$-function. In Section \ref{sec:Uniformity.II} we show that any Euler product, such that the local factors are integrals of the same $V$-function, has rational abscissa of convergence.

\subsection{Acknowledgment} This work is a part of the author's Ph.D. thesis, supervised by Alex Lubotzky. It is a pleasure to thank Alex for introducing the problem to me, for his advices when this work was done, and for his remarks on previous versions of this paper. I have also learned much by talking about this project with Andrei Jaikin, Fritz Grunewald, Michael Larsen, Chris Voll, David Kazhdan, and Udi Hrushovski. I heartily thank them all.

\section{Euler Factorization} \lbl{sec:Euler.Factorization}

\subsection{Notations} In order to remove a layer of unnecessary notational complexity, we assume that the arithmetic lattice $\Ga$ is defined over $\bQ$. That is, we assume that $\Ga =\uG (\bZ _\Sigma)$ where $\Sigma$ is a finite set of prime numbers and $\uG \subset (\underline{GL_n})_{\bZ _\Sigma}$ is a linear algebraic group scheme over $\Spec \bZ _\Sigma$ whose generic fiber is semisimple, simply connected, and connected. The proof for general $\Ga$ is completely analogous.

For every prime $p$ not in $\Sigma$, we denote the group $\uG (\bZ _p)$ by $G_p$. The first congruence subgroup of $G_p$---which we denote by $G_p ^1$---is the kernel of the reduction modulo $p$ homomorphism from $G_p$ to $\uG (\bF _p)$.

In the following, the word `representation' will have several meanings. If the group is discrete, we just mean a complex representation of finite dimension. For profinite groups, a representation should also be continuous (and thus have finite image). If the group is algebraic (or, more generally, pro-algebraic), a representation should be (finite dimensional and) rational. This remark applies also to related notions, such as $\Irr H$ and $\ze _H(s)$.

\subsection{Euler Factorization} In this subsection, we describe without proofs the Euler factorization of $\ze _{\Ga}(s)$, and refer the reader to \cite{LL} for the details. If $\De$ is a finitely generated group, we denote by $\widehat{\De}$ the pro-finite completion of $\De$. For a  discrete group $\De$, the pro-algebraic completion of $\De$ is defined to be a pro-algebraic group $\De^a$, together with a homomorphism $\pi :\De\to \De^a$, such that every finite dimensional representation of $\De$ factors uniquely through $\De^a$. The pro-algebraic completion is unique up to isomorphism, and, by definition, it has the same representations as the group itself. Therefore $\ze _{\De}(s)=\ze _{\De^a}(s)$. Note that in this equality, the left hand side counts all representations of $\De$, whereas the right hand side counts only rational representations of $\De^a$.

Suppose $\Ga =\uG (\bZ _\Sigma)$ satisfies the congruence subgroup property. It was shown in \cite{LL} that there is a finite index subgroup $\De \subset \Ga$ such that the pro-algebraic completion of $\De$ is the direct product of $\widehat{\De }$ and $\uG (\bC )$. Because of the congruence subgroup property, and by making $\De$ smaller if necessary, we may assume that $\widehat{\De}$ is a subgroup (of finite index) of $\prod _{p \not \in \Sigma} \uG (\bZ _p)$. We shall see later (Corollary \ref{cor:finite.index.equal.absc}) that the abscissa of convergence of $\De$ is equal to the abscissa of convergence of $\Ga$. 

Denote the projection from $\widehat{\De}$ to $\uG (\bZ _p)$ by $\pi _p$. Since $\widehat{\De}$ is of finite index in $\prod _{p\not \in \Sigma}\uG(\bZ_p)$, there is a finite set of primes $T$, such that if $p$ does not belong to $T$, then $\pi _p(\widehat{\De})=\uG (\bZ _p)$. We have then
\begin{equation} \lbl{eq:Euler}
\ze _{\De}(s) = \prod _{p\in T} \ze _{\pi _p (\widehat{\De})}(s) \cdot \prod _{p\not \in T} \ze _{G _{p}} (s) \cdot \ze _{\uG (\bC)}(s) .
\end{equation}

We shall call the factor $\ze_{G_p}(s)$ (or $\ze_{\pi _p(\De)}(s)$) the local zeta function at the prime $p$, and call the factor $\ze_{\uG(\bC)}(s)$ the local zeta function at infinity.

The abscissa of convergence for the local zeta function at infinity, $\ze _{\uG (\bC )}(s)$, was computed in \cite{LL}. If $\uG$ has root system $\Phi$ and if we denote the rank of $\Phi$ by $r$, and denote the number of positive roots of $\Phi$ (relative to some choice of simple roots) by $|\Phi _+|$, then the abscissa of convergence of $\ze _{\uG (\bC)}(s)$ is equal to $\frac{r}{|\Phi _+|}$. In particular, the abscissa of convergence of the local zeta function at infinity is rational. 

In \cite{Jai}, the following theorem is proved:

\begin{thm} \lbl{thm:jaikin} For every prime $p$ there are 
\begin{enumerate}
\item A finite set $I_p$.
\item Polynomials $f_i ^p(x) \in \bZ [x]$ with non-negative coefficients for $i\in I_p$.
\item Non-negative integers $n_i$, for $i\in I_p$, and non-negative integers $A_{i,j},B_{i,j}$ for $i\in I_p$ and $1\leq j\leq n_i$
\end{enumerate}
such that
\[
\ze _{G_p}(s)=\sum _{i\in I_p} n_i^{-s} \frac{f_i ^p (p^{-s})}{\prod _j \left (1-p^{-A_{i,j}s+B_{i,j}} \right )}.
\]
The same is true for every finite index subgroup of $G_p$.
\end{thm}

In particular, the abscissa of convergence for every local zeta function is rational. In order to prove that the abscissa of convergence of the `global' zeta function is rational, we shall need to understand the relation between the local zeta functions for different primes. Indeed, this paper is mainly an attempt to give an approximate formula to the local zeta functions, which is uniform in the prime $p$.

\section{Algebraic Preliminaries} \lbl{sec:Algebraic.Preliminaries}

\subsection{Relative Zeta Functions}
Let $H$ be a group and let $K$ be a subgroup of $H$. If $\rho$ is a representation of $K$, we denote its induction to $H$ by $\Ind_K^H\rho$. If $\chi$ is a representation of $H$, we denote its restriction to $K$ by $\Res_K^H\chi$.

\begin{defn} \lbl{defn:rel.zeta} Let $H$ be a group, let $K$ be a normal subgroup of $H$, and let $\tau$ be an irreducible representation of $K$. We denote by $\Irr(H|\tau )$ the set of irreducible representations $\rho$ of $H$ such that $\tau$ is a sub-representation of $\Res _K ^H \rho$ (or equivalently, such that $\rho$ is a sub-representation of $\Ind _K ^H \tau$). Note that if $\rho \in \Irr(H|\tau)$, then $\dim\tau$ divides $\dim\rho$. 

Let $r_n(H|\tau)$ be the number of representations in $\Irr(H|\tau)$ of dimension $n\cdot \dim\tau$. We define the relative zeta function as
\[
\ze _{H|\tau}(s)=\sum _n r_n(H|\tau)\cdot n^{-s} = \sum _{\rho \in \Irr(H|\tau)}\left (\frac{\dim \rho}{\dim \tau} \right )^{-s}.
\]
\end{defn}

\begin{lem} \lbl{lem:rel.zeta} Let $H$ be a group and let $K$ be a normal subgroup of $H$ of finite index. The group $H$ acts on the set $\Irr(K)$ by conjugation. For every $\tau\in\Irr(K)$ we denote the stabilizer of $\tau$ under this action by $\Stab_H\tau$. Then
\[
\ze_H (s) = \sum _{\tau\in\Irr(K)}\frac{1}{[H:\Stab_H\tau]}(\dim \tau)^{-s} \ze_{H|\tau}(s)
\]
\end{lem}

\begin{proof} Let $E$ be the set of pairs $(\tau,\rho )\in \Irr(K) \times \Irr(H)$ such that $\tau$ is a sub-representation of $\Res _K ^H \rho$. Then
\[
\sum _{(\tau,\rho )\in E} \frac{1}{[H:\Stab _H \tau ]}(\dim \rho )^{-s} = \sum _{\tau \in \Irr(K)}\frac{1}{[H:\Stab_H \tau ]}(\dim \tau )^{-s} \ze _{H|\tau}(s).
\]
On the other hand, for every $\rho \in \Irr(H)$, the set of $\tau \in \Irr(K)$ such that $\tau$ is a sub-representation of $\Res_K ^H \rho$ is a single $H$ orbit and so
\[
\sum _{(\tau ,\rho )\in E} \frac{1}{[H:\Stab _H \tau]}(\dim \rho )^{-s} = \sum _{\rho \in \Irr(H)} (\dim \rho )^{-s} \left ( \sum _{\tau | (\tau ,\rho )\in E}\frac{1}{[H:\Stab _H \tau ]} \right ) = \sum _{\rho \in \Irr(H)}(\dim \rho )^{-s}.
\]
\end{proof}

\begin{lem} \lbl{lem:fin.index} Let $K \subset H \subset L$ be groups. Assume that $H$ is of finite index in $L$ and that $K$ is normal in $L$. Let $\tau \in \Irr(K)$. Then for each $N$,
\[
\frac{1}{[L:H]}\left ( r_1(H|\tau )+\ldots +r_{N/[L:H]}(H|\tau) \right ) \leq r_1(L|\tau)+\ldots +r_N(L|\tau) \leq \left ( r_1(H|\tau)+\ldots +r_N(H|\tau) \right ) [L:H]
\]
and for every $s\in \bR$, if one of $\ze_{H|\tau}(s)$ or $\ze_{L|\tau}(s)$ converges, then so does the other, and 
\[
[L:H]^{-1-s}\ze_{H|\tau}(s) \leq \ze _{L|\tau}(s) \leq [L:H]\cdot \ze _{H|\tau}(s).
\]
\end{lem}

\begin{proof} Consider the bipartite graph whose vertices are $\Irr(L|\tau)\sqcup\Irr(H|\tau)$ and there is an edge between $\rho_1 \in \Irr(L|\tau)$ and $\rho_2 \in \Irr(H|\tau)$ if $\rho_2$ is a sub-representation of $\Res^L _H\rho_1$. Note that
\begin{enumerate}
\item Every vertex has positive degree.
\item The degree of every vertex is bounded by $[L:H]$.
\item If $\rho_1\in\Irr(L|\tau)$ and $\rho_2\in\Irr(H|\tau)$ are connected, then $\dim\rho_1\leq \dim\rho_2 \leq [L:H]\cdot \dim\rho_1$.
\end{enumerate}
Let $\Irr(L|\tau)_N\subset\Irr(L|\tau)$ be the set of representations of dimension less than or equal to $N\dim\tau$, and define similarly the set $\Irr(H|\tau)_N$. The set $\Irr(L|\tau)_N$ is contained in the set of neighbors of $\Irr(H|\tau)_N$, so
\[
|\Irr(L|\tau)_N|\leq|\Irr(H|\tau)_N|\cdot[L:H].
\]
Similarly, the set $\Irr(H|\tau)_{N/[L:H]}$ is contained in the the set of neighbors of $\Irr(L|\tau)_N$, so
\[
|\Irr(H|\tau)_{N/[L:H]}|\leq|\Irr(L|\tau)_N|\cdot[L:H].
\]
This proves the first two inequalities. Similar argument shows the other two.
\end{proof}

\begin{cor} \lbl{cor:finite.index.equal.absc} If $H\subset L$ is a subgroup of finite index, then the abscissae of convergence of $\ze_{H}(s)$ and of $\ze_{L}(s)$ are equal.
\end{cor}

\begin{proof} Take $K$ to be the trivial group in Lemma \ref{lem:fin.index}.
\end{proof}

\subsection{Lie Algebras} \lbl{subsec:Lie.algebras}
There are several notions of Lie algebras, exponential functions, and logarithmic functions for pro-$p$ groups and for finite subgroups of $\GL_n(\bF_p)$. We shall give them all here in order to fix notations. We assume in this section that $\uG \subset \underline{\GL_n}_{\bZ_\Sigma}$ is a group scheme over $\Spec \bZ_\Sigma$. Recall that we denote the group $\uG(\bZ_p)$ by $G_p$, and denote its first congruence subgroup by $G_p^1$.

We start with the Lie algebra of $G_p ^1$. There are three definitions for the Lie algebra $\fg^1_p$ of $G_p ^1$. Fortunately, they coincide for almost all primes.

Let $\underline{\fg} \subset \underline{M_n}_{\bZ _\Sigma}$ be the tangent space at the identity, relative to $\Spec\bZ _\Sigma$. For every $p\notin \Sigma$, the set $\underline{\fg}(\bZ _p) \subset M_n(\bZ _p)$ is closed under addition and under taking commutators. The algebraic Lie algebra of $G_p ^1$ is the set
\[
\{ A \in \underline{\fg}(\bZ _p) | A \equiv 0 \quad \textrm{(mod $p$)} \} ,
\]
together with the addition and Lie brackets induced from $M_n(\bZ _p)$. Using the embedding of $G_p$ into $GL_n(\bZ _p)$, the analytic Lie algebra is the set of all matrices of the form
\[
\log (I-g) = (I-g)+\frac{(I-g)^2}{2} +\frac{(I-g)^3}{3}+\dots
\]
where $g\in G_p ^1$. For the analytic Lie algebra, the Lie algebra operations, i.e. addition and Lie brackets, are the usual addition and commutator of matrices. The last definition, due to Lazard (see \cite[Section 4.5]{DdSMS}) is that the Lie algebra, as a set, is just $G_p ^1$, but the addition and brackets need to be redefined. As stated before, those three definitions give isomorphic Lie algebras for almost all primes. We denote the algebraic Lie algebra of the group $G_p^1$ by $\fg_p^1$.

We shall use all three definitions. The algebraic definition implies that there is formula $\phi (x_{i,j})$ in $n^2$ variables, in the language of valued fields (see Section \ref{sec:Definable.Families}), such that for every prime $p$ and every $A\in M_n(\bQ _p)$, we have $A\in \fg^1_p$ if and only if $\phi (A)$ holds. This will enable us to connect the $\fg^1_p$'s for different primes $p$. The analytic definition is useful in order to treat other pro-p subgroups of $G_p$---we shall promptly do this. The Lazard definition is used in \cite{Jai}, to which we shall refer.

We fix $n$ and let
\[
\sU = \{ g\in M_n(\bZ _p) | \lim _{k\to\infty} (g-I)^k =0 \}
\]
and
\[
\sN = \{ A\in M_n(\bZ _p) | \lim _{k\to \infty} A^k =0 \}
\]
be the sets of pro-unipotent and pro-nilpotent elements respectively. For $g\in \sU$ define $\log (g)$ as the series
\[
\log (g) = (g-I) + \frac{(g-I)^2}{2} + \frac{(g-I)^3}{3} + \dots .
\]
For $A\in \sN$ define $\exp (A)$ as the series
\[
\exp (A) = I+A+\frac{A^2}{2!} + \frac{A^3}{3!} + \dots .
\]

\begin{lem} If $p>2n$ then the series defining $\log$ and $\exp$ converge, and the functions $\log, \exp$ are inverses. Moreover, if $A,B\in \sN$, and the reductions mod $p$, $\overline{A},\overline{B}$, are in a nilpotent Lie subalgebra of $M_n(\bF_p)$, then the Campbell Hausdorff formula holds:
\begin{equation} \lbl{eq:CH}
\log (\exp (A)\cdot \exp (B)) = \sum _{m=1} ^{\infty} \frac{(-1)^m}{m}\sum _{r_i+s_i >0}\frac{\left ( \sum _{i=1}^m (r_i+s_i) \right ) ^{-1}}{r_1!\cdot s_1! \cdot \ldots \cdot r_m!\cdot s_m!} R_{r_1,s_1,\ldots ,r_m,s_m}(A,B)
\end{equation}
where $R_{r_i,s_i} (A,B)$ is defined by
\[
R_{r_1,s_1,\ldots ,r_m,1}(A,B)=(ad(A))^{r_1}(ad(B))^{s_1} \cdots (ad(A))^{r_m}(B) ,
\]
\[
R_{r_1,s_1,\ldots ,1,0}(A,B)=(ad(A))^{r_1}(ad(B))^{s_1} \cdots (ad(B))^{r_{m-1}}(A) ,
\]
and $R_{r_i,s_i} (A,B)=0$ otherwise.
\end{lem}

\begin{proof} If $A\in \sN$ then $A^n$ is divisible by $p$. Therefore for every $N$, $A^N$ is divisible by $p^{\lfloor \frac{N}{n} \rfloor}$. As the maximal power of $p$ that divides $N!$ is $\lfloor \frac{N}{p} \rfloor + \lfloor \frac{N}{p^2} \rfloor + \ldots < \frac{2N}{p}$, we get that if $p>2n$, then the term $\frac{A^N}{N!}$ is divisible by $p^{\lfloor \frac{N}{n} \rfloor -\frac{2N}{p}}$. We get that $v_p \left ( \frac{A^N}{N!} \right )$ tends to infinity as $N$ tends to infinity. Therefore, the series defining $\exp$ is convergent. The same argument shows that the series defining $\log$ is convergent.

Similarly, if $\overline{A},\overline{B}$ are contained in a nilpotent Lie subalgebra of $\fg \fl _n(\bF _p)$, then $R_{r_1,s_1,\ldots ,r_m,s_m}(A,B)$ is divisible by 
\[
p^{\frac{r_1+s_1+\ldots +r_m+s_m}{n}} ,
\]
whereas the maximal power of $p$ that divides $r_1!s_1!\ldots r_m!s_m!$ is less than
\[
\frac{2(r_1+s_1+\ldots +r_m+s_m)}{p}.
\]
So the right hand side of ($\ref{eq:CH}$) is convergent, and therefore is equal to the left hand side.
\end{proof}

\begin{defn} Let $R\subset G_p$ be a pro-$p$ subgroup of $G_p$ such that $G_p^1\subset R$. Since $\overline{R}$---the reduction of $R$ modulo $p$---is a $p$-subgroup of $GL_n(\bF_p)$, we know that every element in $\overline{R}$ is unipotent. Hence $R\subset\sU$. We define the Lie algebra of $R$ to be the image of $R$ under the map $\log$, and denote it by $\Lie(R)$.
\end{defn}
Note that $\overline{\Lie(R)}\subset\bM_n(\bF_p)$ is a nilpotent Lie algebra.

There is yet another notion of Lie algebras, this time for subgroups of $\GL_n(\bF_p)$. It is taken from \cite{No}. Assume $p>2n$ and let $\Upsilon\subset\GL_n(\bF_p)$. If $\ga\in\Upsilon$ is an element of order $p$, then $(\ga-I)^n=0$. We define
\[
\log(\ga)=(\ga-I)+\frac{(\ga-I)^2}{2}+\ldots+\frac{(\ga-I)^{n-1}}{n-1}.
\]
The Lie algebra of $\Upsilon$ is the set
\[
\Lie(\Upsilon)=\bF_p-span\{\log(\ga) | \ga\in\Upsilon, \textrm{ the order of $\ga$ is $p$} \}
\]
The set $\Lie(\Upsilon)$ is shown in \cite{No} to be closed under commutators.

\subsection{Orbit Method}
Recall that for a locally compact abelian group $A$, the Pontjagin dual of $A$---which we denote by $A^\vee$---is the set of all continuous homomorphisms from $A$ to the circle group $S^1=\{z\in\bC\qq|\qq|z|=1\}$.

Let $R\subset G_p$ be a pro-$p$ subgroup such that $G_p^1\subset R$. Let $\fr=\Lie(R)$. The group $R$ acts on the (additive) group $\fr$ by conjugation, and therefore acts on the Ponrjagin dual $\fr ^{\vee}$. We call this action the {\em coadjoint} action and denote it by $Ad^*$. Concretely, it is given by
\[
(Ad^*(g)\theta)(X)=\theta (X^{g^{-1}}) = \theta (gXg^{-1}).
\]

\begin{thm} \lbl{thm:orbit.method} Given $\uG$, there is an integer $p_0$, such that if $p>p_0$ is a prime, and if $Q\subset R\subset G_p$ are pro-$p$ subgroups of $G_p$ with Lie algebras $\fq \subset \fr$ respectively, then the following hold
\begin{enumerate}

\item \lbl{cl:orbit.1} There is a bijection $\Xi_R$ between $Ad^*(R)$ orbits on $\fr ^{\vee}$ and irreducible representations of $R$. If $\theta \in \fr ^{\vee}$, we shall write $\Xi_R(\theta )$ instead of $\Xi_R (Ad^* (R)\theta )$.

\item \lbl{cl:orbit.2} The character of $\Xi_R(\theta )$ is given by
\[
\chi _{\Xi_R(\theta)}(g) = \frac{1}{|Ad^*(R)\theta|^{1/2}} \sum _{\phi \in Ad^*(R)\theta} \phi (\log (g))
\]

\item \lbl{cl:orbit.3} If $\theta \in \fr ^{\vee}$, then the dimension of $\Xi_R(\theta)$ is $|Ad^*(R)\theta |^{1/2}$.

\item \lbl{cl:orbit.4} If $\theta \in \fr ^{\vee}$, and $\tau \in \fq ^{\vee}$, then $\Xi_Q(\tau )$ is a sub-representation of $\Res_S ^R \Xi_R (\theta )$ if and only if there is $g\in R$ such that $\tau = Ad^* (g)\theta |_{\fq}$.

\end{enumerate}
\end{thm}

\begin{proof} The proof of (\ref{cl:orbit.1}) and (\ref{cl:orbit.2}) is identical to the proof of theorem 1.1 in \cite{Ho}, using the fact that $\fr$ is closed under addition and brackets, and using the Campbell Hausdorff formula.

(\ref{cl:orbit.3}) follows from (\ref{cl:orbit.2}) by evaluating the character at $1$.

(\ref{cl:orbit.4}): By (\ref{cl:orbit.2}), for every $g\in Q$ the evaluation at $g$ of the characters of $\Xi_Q(\tau)$ and $\Res_Q ^R \Xi_R(\theta)$ are
\[
\chi _{\Xi_Q(\tau)}(g) = \frac{1}{|Ad^*(Q)\tau|^{1/2}} \sum _{\phi \in Ad^*(Q)\tau} \phi (\log (g))
\]
and
\[
\chi _{\Res_Q ^R \Xi_R(\theta)} = \frac{1}{|Ad^*(R)\theta|^{1/2}} \sum _{\psi \in Ad^*(R)\theta} \psi (\log (g)) .
\]
The map $\exp :\fq \to Q$ is a measure preserving bijection and hence
\[
(\Xi_Q(\tau),\Res_Q ^R \Xi_R(\theta)) = \int _Q \chi _{\Xi_Q(\tau)} (g) \cdot \overline{\chi _{\Res_Q ^R \Xi_R(\theta)}(g)} dg = 
\]
\[
=\frac{1}{|Ad^*(Q)\tau|^{1/2}} \cdot \frac{1}{|Ad^*(R)\theta|^{1/2}} \sum _{\phi \in Ad^*(Q)\tau} \sum _{\psi \in Ad^*(R)\theta} \int _{\fq} \phi (X)\cdot \overline{\psi |_{\fq}(X)}dX .
\]
Every $\phi$ and $\psi |_{\fq}$ in the above sum are one dimensional characters of $\fq$, and by orthogonality of characters, 
\[
\int _\fq \phi (X) \overline{\psi |_\fq (X)}dX = \left \{ \begin{matrix} 1 & \textrm{if $\phi=\psi |_\fq$} \\ 0 & \textrm{if $\phi \neq \psi |_\fq$} \end{matrix} \right . .
\]
The claim follows immediately from this.
\end{proof}

\subsection{Subgroups of $GL_n(\bF_p)$}
We review some definitions from \cite{No}, and advise the reader to have a copy in hand. We fix a natural number $n$ and a prime $p$. Let $\Upsilon$ be a subgroup of $\GL_n(\bF_p)$. We shall denote by $\Upsilon^+$ the subgroup of $\Upsilon$ which is generated by the $p$-elements of $\Upsilon$.

If $L\subset\M_n(\bF_p)$ is a Lie subalgebra, we denote by $\underline{\exp L}$ the algebraic group generated by the one parameter subgroups
\[
t\mapsto \exp(tX)
\]
for all nilpotent $X\in L$.

If $\Upsilon$ is a subgroup of $\GL_n(\bF_p)$, we define 
\[
\widetilde{\Upsilon}=\underline{\exp(\Lie(\Upsilon))}
\]
where $\Lie(\Upsilon)$ was defined in Subsection \ref{subsec:Lie.algebras}.

For a subset $S\subset\Upsilon$, we denote the subgroup generated by $S$ by $\langle S\rangle$. For an algebraic group $G$ we denote $Lie(G)$ the Lie algebra of $G$. If $G$ is defined over $\bF_p$, then $Lie(G)$ can be thought of as a subalgebra of $\fg\fl_n(\bF_p)$. If $L\subset\fg\fl_n(\bF_p)$ we denote by $\exp L$ the set of elements of the form $\exp X$ for $X\in L$ a nilpotent element.

\begin{prop} \lbl{prop:bounded.generation.algebraic} For every $n$ there is an $N$ such that if $L\subset gl _n(\bF _p)$ is Lie algebra that is generated by nilpotents, then $\langle \exp L \rangle =(\exp L)^N$. Moreover, there are elements $X_1,\ldots ,X_N \in L\cap N_n(\bF _p)$ such that
\[
\langle \exp L \rangle = \langle \exp X_1 \rangle \cdot \ldots \cdot \langle \exp X_N \rangle
\]
\end{prop}

\begin{proof} The claims in the proposition are trivial if $p$ is bounded. For the proof, we shall assume that $p$ is large enough, and so we can use the results of \cite{No}. Also, the first statement clearly follows from the second, so we prove the second claim. 

Let $R\subset L$ be the unipotent radical of $L$. We define algebraic groups
\[
A=\underline{\exp L} \quad \textrm{and} \quad B=\underline{\exp R} .
\]
The algebraic group $B$ is normal in $A$, and therefore the group $B (\bF _p)$ is normal in $A(\bF _p)$. We have an exact sequence
\[
0\to B(\bF _p) \to A(\bF _p) \to (A/B)(\bF _p) \to H^1(\bF _p,B).
\]
Since $B$ is unipotent, the Galois cohomology group, $H^1(\bF _p,B)$, vanishes, and hence $A(\bF _p)/B(\bF _p)=(A/B)(\bF _p)$. Since $B(\bF_p)$ is a $p$-group we get that $B(\bF_p) \triangleleft A(\bF_p) ^+$, and that $A(\bF_p)^+ /B(\bF_p)=(A/B)(\bF_p)^+$.

The first part of thorem A of \cite{No} implies that $Lie(A)=L$ and $Lie(B)=R$. Hence $Lie(A/B)=L/R$. Since $L$ is nilpotently generated, so is $L/R$. Hence the Lie algebra $L/R$ is semisimple (and not only reductive). In this case, there are nilpotent elements $\overline{X_1},\ldots ,\overline{X_M}\in Lie(A/B)$ (where $M$ depends only on $n$), such that
\begin{equation} \lbl{eq:b.g.s.s}
(A/B)(\bF _p)^+=\langle \exp \overline{X_1}\rangle \cdot \ldots \cdot \langle \exp \overline{X_M} \rangle .
\end{equation}
By induction on the nilpotency class of $R$, there are nilpotent elements $Y_1,\ldots ,Y_K \in R$ such that
\begin{equation} \lbl{eq:b.g.unip}
B(\bF _p)=\langle \exp Y_1 \rangle \cdot \ldots \cdot \langle \exp Y_K \rangle .
\end{equation}
Choose $X_i \in L$ such that $X_i+R=\overline{X_i}$. From (\ref{eq:b.g.s.s}) and (\ref{eq:b.g.unip}) we get that 
\[
A(\bF_p)^+ = \langle \exp {X_1}\rangle \cdot \ldots \cdot \langle \exp {X_M} \rangle \cdot \langle \exp Y_1 \rangle \cdot \ldots \cdot \langle \exp Y_K \rangle .
\]
It remains to show that $A(\bF_p)^+=\langle \exp L \rangle$. Clearly, if $X\in L$ then $\exp X\in A (\bF_p)$, and so $\langle \exp L \rangle \subset A(\bF_p)^+$. For the converse, suppose that $u\in A(\bF_p)$ is an element of order $p$. Denoting $X=\log u$, we get
\[
X=\log u \in \langle \log A(\bF_p) \rangle \stackrel{(1)}{=} Lie \left ( \widetilde{A(\bF_p)} \right ) \stackrel{(2)}{=} Lie (A) = Lie (\underline{\exp L}) \stackrel{(3)}{=} L,
\]
where $(1)$ follow from first part of theorem B of \cite{No}, $(2)$ follows from the second part of the same theorem, and $(3)$ follows from the first part of theorem A of \cite{No}.

Since $A(\bF_p)^+$ is generated by the $p$-elements in $A(\bF_p)$, we get that $A(\bF_p)^+\subset\langle\exp L\rangle$.
\end{proof}

We have the following corollary, which will not be used in the rest of the article:

\begin{cor} \lbl{cor:b.g.general} For every $n$ there is $N$ such that for every prime $p$ and a group $G\subset GL_n(\bF _p)$, there are elements $x_1,\ldots,x_N \in G$ such that
\[
G=\langle x_1 \rangle \cdot \ldots \cdot \langle x_N \rangle .
\]
\end{cor}

\begin{proof} By theorem B of \cite{No}, $G^+=\widetilde{G}(\bF_p)^+$. Every nilpotent element in $Lie (\widetilde{G})$ is of the form $\log g$, where $g\in \widetilde{G}(\bF _p)^+=G^+$. Hence, by Proposition \ref{prop:bounded.generation.algebraic}, there are elements $g_1,\ldots ,x_M\in G^+$ such that 
\[
G^+=\langle g_1 \rangle \cdot \ldots \cdot \langle g_M \rangle .
\]
By theorem C of \cite{No}, there is an abelian group $H\subset G$ such that $HG^+$ is normal and of bounded index in $G$. Moreover, by the proof of the theorem, there is a lifting of $H\hookrightarrow GL_n(\bF _p)$ to $H\hookrightarrow GL_n(\bZ _p)$. Hence $H$ is a finite abelian subgroup of $GL_n(\bQ_p)$, and hence is a product of at most $n$ cyclic groups. It follows that one can find $h_1,\ldots ,h_n\in H$ such that
\[
H=\langle h_1 \rangle \cdot \ldots \cdot \langle h_n \rangle .
\]
Finally, there are elements $z_1,\ldots,z_{\log _2 [G:HG^+]}\in G$ such that
\[
G/HG^+ = \langle z_1HG^+ \rangle \cdot \ldots \cdot \langle z_{\log _2 [G:HG^+]}HG^+ \rangle .
\]
Putting it together,
\[
G= \langle g_1 \rangle \cdot \ldots \cdot \langle g_M \rangle \cdot \langle h_1 \rangle \cdot \ldots \cdot \langle h_n \rangle \cdot \langle z_1 \rangle \cdot \ldots \cdot \langle z_{\log _2 [G:HG^+]} \rangle .
\]
\end{proof}

\subsection{Extensions of Representations} \lbl{subsec:ext.of.reps}
In general, if $S$ is a group, $V\triangleleft S$ is a normal subgroup, and $\rho$ is a representation of $V$, then the relative zeta function $\ze_{S|\rho}(s)$ is different from the representation zeta function of the quotient, $\ze_{S/V}(s)$. There is, however, one important case in which they are equal:

\begin{defn} Let $S$ be a group and $V\triangleleft S$ be a normal subgroup. Let $\rho$ be a representation of $H$. We say that $\rho$ is extendible to $S$, if there is a representation $\chi$ of $S$ such that $\Res_V^S\chi=\rho$. The representation $\chi$ is called an extension of $\rho$ to $S$.
\end{defn}

Suppose $\rho\in\Irr V$ is extendible to $S$, and let $\chi$ be an extension of $\rho$ to $S$. Every representation $\tau$ of $S/V$ can be thought of as a representation of $S$ by composition with the quotient map $S\to S/V$. We have a map
\begin{equation} \lbl{eq:rel.Irr}
\Irr(S/V)\longrightarrow\Irr(S|\rho)
\end{equation}
\[
\tau\mapsto\tau\otimes\chi
\]

\begin{prop} \lbl{prop:Isaccs.6.16} (\cite[Theorem 6.16]{Is}) If $V\triangleleft S$, $\rho\in\Irr V$, and $\chi\in\Irr S$ is an extension of $\rho$ to $V$, then the map (\ref{eq:rel.Irr}) is a bijection. Therefore $\ze_{V|\rho}(s)=\ze_{S/V}(s)$.
\end{prop}

Extensions of representations are tightly connected to the second cohomology group of the quotient. The setting is as follows: We have a group $S$, a normal subgroup $V\triangleleft S$, and an irreducible representation $\rho$ of $V$. By Clifford's theory, a necessary condition for the extendability of $\rho$ is that $S$ fixes the representation $\rho$\footnote{I.e. that for every $g\in S$, the representation $\rho^g$, defined by $v\mapsto \rho(g^{-1}vg)$, is equivalent to $\rho$.}. Assuming this, we construct an element in the second cohomology group $H^2(S/V,\bC^\times)$.

Let $M$ be a $V$-module that gives rise to the representation $\rho$. Choose a transversal $T$ to $V$ inside $S$, such that $1\in T$. For every $t\in T$, the $V$-modules $M$ and $tM$ are isomorphic. We choose an isomorphism $P_t:M\to tM$, and for $t=1$ we put $P_1=\Id$. Every element of $S$ can be written as $tv$, where $t\in T$ and $v\in V$. We define $P_{tv}:M\to tM$ as $P_{tv}(m)=P_t(v\cdot m)$. It can be easily checked that for any $g_1,g_2\in S$, the operator 
\[
P_{(g_1g_2)}^{-1} \circ P_{g_1}\circ P_{g_2}:M \to M
\]
is a morphism of $V$ modules, and hence it is a multiplication by a scalar, which we denote by $\al(g_1,g_2)$. Note that the value of $\al(g_1,g_2)$ depends only on the cosets $g_1V,g_2V$. The function $\al$ is a 2-cocycle, and we denote its image in the second cohomology of $S/V$ by $\be$. Although the cocycle $\al$ depends on the choices of $T$ and $P_t$, the cohomology class $\be$ does not. By \cite[Theorem 11.7]{Is}, the representation $\rho$ is extendible to $S$ if and only if $\be$ is trivial.

We will be interested in the case that $V$ is a pro-$p$ group. In this case, we have
\begin{prop} \lbl{prop:extension.p.group} Let $S$ be a profinite group, let $V\triangleleft S$ be a normal pro-$p$ subgroup of finite index, and let $\rho$ be an irreducible representation of $V$. If $\be\in H^2(S/V,\bC^\times)$ is the cohomology class attached to $S,V,\rho$, then $\be$ is a $p$-element in $H^2(S/V,\bC^\times)$.
\end{prop}

\begin{proof} Fix volume forms on the $tM$'s. Let $t_1,t_2\in T$ and suppose $v\in V$ is such that $t_1t_2v\in T$. By taking determinants we get
\[
\al(t_1,t_2)^{\dim(M)} = \det(P_{t_1})\det(P_{t_2})\det(P_{t_1t_2v})^{-1} \det \rho(v)^{-1}.
\]
We can choose the $P_t$'s for $t\in T$ to have determinant 1. Since $\rho$ is an irreducible representation of a pro-$p$ group, $\dim(M)$ is a power of $p$, and $\det \rho(v)$ is a $p^n$-root of unity for some $n$. Therefore $\al(t_1,t_2)$ is a $p^m$-root of unity. It follows that the order of $\be$ is a power of $p$. 
\end{proof}

The cohomology groups of finite quasi-simple groups are well known. In particular, we have

\begin{prop} \lbl{prop:cohomology.finite.simple} For every $r\in\bN$ there is $c(r)\in\bN$ such that if $\Theta$ is a quasi-simple group of Lie rank $r$, then the order of the group $H^2(\Theta,\bC^\times)$ is less than $c(r)$. The same is true for the first cohomology groups $H^1(\Theta,\bC^\times)$.
\end{prop}

We shall use Propositions \ref{prop:extension.p.group} and \ref{prop:cohomology.finite.simple} for extensions of $p$ groups by finite quasi-simple groups, where the rank of the finite quasi-simple group is bounded, and $p$ is large. In this case there are no $p$-elements in the second cohomology group, and therefore the relative zeta function is equal to the representation zeta function of the finite quasi-simple group. 

\subsection{Zeta Functions of Finite Reductive Groups}
The representation zeta functions of the finite simple groups of Lie type were studied in \cite{LiSh} using the Deligne-Lusztig theory. 

Let $G$ be a connected, simply connected, and simple algebraic group defined over $\bF_p$. Let $T\subset G$ be a maximal torus defined over $\bF_p$. Choose a Borel subgroup $B\subset G$, not necessarily defined over $\bF_p$, and let $U$ be the unipotent radical of $B$.

Let $F:G\to G$ be the Frobenius map. Recall that the Lang map $L:G\to G$ is the map
\[
g\mapsto g^{-1}\cdot F(G).
\]

The group $G(\bF_p)\times T(\bF_p)$ acts on the variety $L^{-1}(U)$ by $(g,t)(x)=gxt^{-1}$. Therefore, for each $i$, the $i$-th \'etale cohomology with compact support, $H^i_c(L^{-1};\bC)$\footnote{To be more precise, for every prime $\ell\neq p$, we have the cohomology groups $H^i_c(L^{-1};\overline{\bQ_\ell})$. But $\overline{\bQ_\ell}$ is isomorphic to $\bC$.}, is a $G(\bF_p)\times T(\bF_p)$-bimodule.

\begin{defn} The Deligne--Lusztig induction of a character $\theta$ of $T(\bF_p)$ is the $\theta$-isotypic component in the virtual\footnote{Virtual means that we are taking formal linear combinations of representations. The result lives in the K group of the category of representations.} module
\[
\sum (-1)^i H^i_c(L^{-1}(U);\bC).
\]
This is a virtual representation of $G(\bF_p)$, it is independent of the choice of $B$, and we denote it by $R_T^G\theta$.
\end{defn}

\begin{lem} For fixed $T$ and $\theta$, all irreducible components of $R_T^G\theta$ have the same central character.
\end{lem}

\begin{proof} For a variety $V$ and $f\in\Aut(V)$, let
\[
\cL(g,V)=\sum (-1)^i \trace(f|H^i_c(V;\bC)).
\]
By the definition, the character of $R_T^G\theta$ is
\[
\trace(R_T^G\theta(g))=\frac{1}{|T(\bF_p)|}\sum _{t\in T(\bF_p)}\cL((g,t),L^{-1}(U))\cdot \theta(t)^{-1}.
\]

If $z\in Z(G(\bF_p))$ and $t\in T(\bF_p)$, then the order of $(z,t)$ as an automorphism of $L^{-1}(U)$ is prime to $p$. Therefore, by \cite[Proposition 10.14]{DM}, $\cL((z,t),L^{-1}(U))=\cL((1,1),(L^{-1}(U))^{(z,t)})$. If $z\neq t$, then $L^{-1}(U)^{(z,t)}=\emptyset$ and $\cL((z,t),L^{-1}(U))=0$. If $z=t$, then $L^{-1}(U)^{(z,t)}=L^{-1}(U)$. Therefore we get
\[
\trace(R_T^G\theta(z))=\frac{1}{|T(\bF_p)|}\cL((1,1),L^{-1}(U))\cdot \theta(z)^{-1}.
\]
Since
\[
\trace(R_T^G\theta(1))=\frac{1}{|T(\bF_p)|}\cL((1,1),L^{-1}(U)),
\]
the lemma follows.
\end{proof}

The following is a slight generalization of \cite[Theorem 1.7]{LiSh}, which we will need in the following:

\begin{lem} \lbl{lem:simple.zeta.function} Let $G$ be a simple group scheme. There is a natural number $N$ such that for every $0\leq a<N$ we have
\begin{enumerate}
\item The isomorphism type of the center of $G(\bF_p)$ is the same, for almost all primes congruent to $a$ modulo $N$. Denote this group by $A_a$.
\item For every $\om\in A_a^\vee$, there are polynomials $P_1,\ldots,P_N,Q_1,\ldots,Q_N$ such that for almost all primes $p$ that are congruent to $a$ modulo $N$ we have
\[
\ze_{G(\bF_p)|\om}(s)=\sum P_i(p)\cdot (Q_i(p))^{-s}.
\]
\end{enumerate}
\end{lem}

\begin{proof} (1) is well known. Let $G^*$ be the dual algebraic group to $G$, as defined in \cite[Definition 13.10]{DM}. The representations of $G(\bF_p)$ are partitioned into Lusztig cells, $\cE(G(\bF_p),(s))$, indexed by semi-simple conjugacy classes $(s)\subset G^*(\bF_p)$. Each cell consists of the irreducible components of the representation $R_T^G\theta$, where the pair $(T,\theta)$, consisting of a maximal torus $T\subset G$ defined over $\bF_p$ and a character $\theta$ of $T$, is attached to the conjugacy class $(s)$ by \cite[Proposition 13.13]{DM}. By \cite[Theorem 13.23 and Remark 13.24]{DM}, for every $s$, there is a bijection $\psi_s$ between $\cE(G(\bF_p),(s))$ and $\cE(C_{G*(\bF_p)}(s),1)$ such that
\[
\dim\rho=\frac{|G(\bF_p)|_{p'}}{|C_{G*}(s)(\bF_p)|_{p'}}\dim\psi_s(\rho),
\]
where $|X|_{p'}$ denotes the largest integer prime to $p$ that divides $|X|$.

If $s\in G^*(\bF_p)$ is semi-simple, then $C_{G^*}(s)$ is a reductive subgroup of $G^*$ of maximal rank and $C_{G^*(\bF_p)}(s)=C_{G^*}(s)(\bF_p)$. There are finitely many subgroup schemes $C_1,\ldots,C_K\subset G^*$ such that for any prime $p$ and every semi-simple $s\in G^*(\bF_p)$, $C_{G^*}(s)$ is conjugate to one of the $C_i$'s. Moreover, for every $i$ there is a polynomial $F^1_i(x)\in\bQ[x]$ such that for every $p$ we have 
\[
F^1_i(p)=\frac{|G(\bF_p)|_{p'}}{|C_i(\bF_p)|_{p'}}.
\]

By looking at the table of unipotent characters, we see that there are polynomials $F^2_{i,j}(x)\in\bQ[x]$, such that for every $p$, the degrees of the unipotent representations of $C_i(\bF_p)$ are $F^2_{i,1}(p),\ldots,F^2_{i,M}(p)$.

Finally, by a similar argument to \cite[Lemma 4.3]{LiSh}, the number of conjugacy classes $(s)\subset G^*(\bF_p)$ such that \begin{enumerate}
\item $C_{G^*}(s)$ is conjugate to $C_i(\bF_p)$.
\item $\theta|_{A}=\om$, where $(T,\theta)$ is the pair associated to $(s)$.
\end{enumerate}
is of the form $F^3_{i,\om,p}(p)$, where $F^3_{i,\om,p}(x)\in\bQ[x]$ depends only on $i,\om$, and the residue class of $p$ modulo some fixed integer $N$. We get that
\[
\ze_{G(\bF_p)|\om}(s)=\sum _{i=1}^K F^3_{i,\om,p}(p)\cdot F^1_i(p)\cdot \sum _{j=1}^M F^2_{i,j}(p)^{-s}.
\]
\end{proof}

\begin{prop} \lbl{prop:zeta.semisimple.finite} Let $G$ be a semisimple algebraic group scheme over $\bZ_S$. There is a natural number $N$, and for each $0\leq a<N$, there are two sequences of polynomials 
\[
P_1(x),\ldots,P_{k_a},Q_1(x),\ldots,Q_{k_a}(x)\in\bQ[x]
\]
such that for every prime $p$, which is congruent to $a$ modulo $N$ and not in $S$, we have
\[
\ze_{G(\bF_p)}(s)=\sum_{i=1}^{k_a}P_i(p)\cdot Q_i(p)^{-s}.
\]
\end{prop}

\begin{proof} Let $G$ be a semisimple algebraic group scheme over $S$. There are simple algebraic group schemes $G_1,\ldots,G_n$ such that for every $p$ we have a central extension
\[
1\to Z_p \to \prod G_i(\bF_p) \to G(\bF_p) \to 1.
\]
Moreover, the isomorphism classes of $Z_p$ and of the centers of $G_i(\bF_p)$ are constant if we fix the residue class of $p$  modulo some $N$. Fix such a residue class and let $\Om$ be the collection of tuples $(\om_1,\ldots,\om_n)$ such that
\begin{enumerate}
\item $\om_i$ is a character of the center of $G_i(\bF_p)$.
\item $\om_1\times\ldots\times\om_n$ is trivial on $Z_p$.
\end{enumerate}
We have that
\[
\ze_{G(\bF_p)}(s)=\sum_{(\om_1,\ldots,\om_n)\in\Om}\ze_{G_1(\bF_p)|\om_1}(s)\cdot\ldots\cdot\ze_{G_n(\bF_p)|\om_n}(s).
\]
By Lemma \ref{lem:simple.zeta.function}, the proposition is proved.
\end{proof}

\subsection{Equivalence of Euler Products}
\begin{defn} \lbl{defn:equiv} Let $(\ze_n(s))_n$ and $(\xi_n(s))_n$ be two sequences of Dirichlet series with non-negative coefficients. We say that the sequences $(\ze_n(s))_n$ and $(\xi_n(s))_n$ are {\em equivalent}, and we write $(\ze _n(s))_n \sim (\xi _n(s))_n$, if there is a constant $C>0$ such that for every $n$ and every $s$, which is greater than the abscissae of convergence of all $\xi _n(s),\ze _n(s)$,
\[
C ^{-1-s} \xi _n (s) \leq \ze _n (s) \leq C ^{1+s}\xi _n (s).
\]
\end{defn}

\begin{lem} \lbl{lem:equiv.Euler} Suppose $(\ze _n(s))_n$ and $(\xi _n(s))_n$ are two sequences of Dirichlet series with non-negative coefficients and constant terms equal to zero, and suppose that $(\ze _n (s))_n \sim (\xi _n(s))_n$. Then the abscissae of convergence of the products
\[
\prod _n (1+\ze _n (s)) \qq \textrm{and} \qq \prod _n (1+\xi _n (s))
\]
are equal.
\end{lem}

\begin{proof} Suppose $s$ is greater than the abscissa of convergence of $\prod _n (1+\ze _n (s))$. Then for every $n$, $s$ is greater than the abscissa of convergence of $\ze _n$, and the sum $\sum _n \ze _n (s)$ converges. By the assumption, $s$ is greater than the abscissa of convergence of $\xi _n$ for every $n$, and the sum $\sum \xi _n (s)$ converges. Therefore, $s$ is greater than the abscissa of convergence of $\prod (1+\xi _n(s))$. By symmetry, the abscissae of convergence of $\prod (1+\ze _n(s))$ and $\prod (1+\xi _n(s))$ are equal.
\end{proof}

\subsection{Resolution of Singularities} \lbl{subsec:resolution.of.singularities}
In this subsection we remind the notions of resolution of singularities and reduction modulo $p$ of a scheme defined over the rationals.

We start with the notion of (embedded) resolution of singularities. We shall work over the field $\bQ$ of rational numbers. Given a polynomial $P(x)\in \bQ [x_1 ,\dots ,x_n]$, an embedded resolution of $P(x)$ is a pair $(Y_\bQ,h)$ where $Y_\bQ$ is a smooth subvariety of $\bP_{\bA_\bQ^n}^k$, and $h$ is the restriction of the natural projection $\bP _{\bA_\bQ^n}^k\to\bA _\bQ^n$ to $Y_\bQ$, such that if we denote by $D$ the subscheme defined by $P(x)$, then
\begin{enumerate}
\item The restriction of $h$ to $Y_\bQ\setminus h^{-1}(D)$ is an isomorphism onto $\bA^n_\bQ\setminus D$.
\item $h^{-1}(D)$ is a divisor with normal crossings.
\end{enumerate}
By a well known theorem of Hironaka, every polynomial (over a field of characteristics 0) has a resolution of singularities.

The second notion we wish to remind the reader is that of reduction mod $p$ of a variety. Let $Y_\bQ\subset\bP_{\bA _\bQ^n}^k$ be a variety. Consider $\bP _{\bA _{\bQ}^n}^k$ as an open subset of $\bP_{\bA _{\bZ}^n}^k$. Define $Y_\bZ$ to be the scheme theoretic closure of $Y_\bQ$ inside $\bP_{\bA_\bZ^n}^k$. The reduction mod $p$ of $Y_\bQ$ is the fiber product $Y_{\bZ}\times _{\Spec(\bZ)}\Spec(\bF_p)$ \footnote{A more elementary description, which is true for almost all primes is the following: Suppose $Y_\bQ$ is defined by the polynomial equations $Q_1(x)=\dots =Q_m(x)=0$ where $Q_i(x)$ are polynomials with rational coefficients. For almost all primes $p$, the denominators of the coefficients of $Q_i(x)$ are not divisible by $p$ and so we can consider the reduction $\overline{Q_i}(x)$ of $Q_i(x)$ mod $p$. Then $Y_{\bF _p}$ is the variety defined by the equations $\overline{Q_1}(x)=\dots \overline{Q_m}(x)=0$.}. 

\begin{defn} Let $P(x)\in\bQ[x_1,\ldots,x_n]$ be a polynomial. Let $(Y_\bQ,h)$ be a resolution of singularities of $P(x)$. We denote the irreducible components of $(h^{-1}(D))_{red}$ be $E_1,\ldots,E_m$. We say that $(Y,h)$ has good reduction modulo $p$ if the following conditions hold:
\begin{enumerate}
\item $Y_{\bF _q}$ is smooth.
\item $\overline{E_i}$ are smooth, and $\cup \overline{E_i}$ has normal crossings.
\item $\overline{E_i}$ and $\overline{E_j}$ do not have a common irreducible component if $i\neq j$.
\end{enumerate}
\end{defn}

It is easy to see that if $(Y_\bQ,h)$ is a resolution of singularities then this resolution has a good reduction modulo almost all primes.

\section{Definable Families} \lbl{sec:Definable.Families}
\subsection{Definable Sets in $\cT _f$} \lbl{subsec:def.set} We shall work with several different logical theories (and languages). Recall that the language of rings, $\cL _{Rings}$, is the first order language which has constant symbols 0,1, has only equality as a relation, and has two function symbols: addition and multiplication. We let the theory $\cT _f$ (the theory of fields) be the collection of all sentences in $\cL _{Rings}$ that hold for all fields.

It will also be useful to work over different bases. If $R$ is an integral domain, we denote by $\cL _{Rings}(R)$ the language $\cL _{Rings}$ together with constant symbols for the elements of $R$. The theory $\cT _f (R)$ consists of all sentences of $\cL _{Rings}(R)$ that hold for all fields containing $R$. In particular, it contains all relations that hold between the elements of $R$.

By a {\bf $\cT _f$-definable set} we shall mean a formula $\phi (x)$ in the language $\cL _{Rings}$ (here and in the following we shall use $x$ to denote a tuple of variables of unspecified length). Let $X$ be a definable set that corresponds to the formula $\phi (x)$. Given a model $L$ of $\cT _{f}$ (i.e. a field) we define the set of $L$-solutions of $X$ as
\[
X(L)=\phi (L) :=\{ a\in L^n \qq | \qq \phi (a)\}.
\]

Examples of definable sets are the affine space $\bA ^n$ defined by the formula $\phi (x_1 ,\dots ,x_n):=`0=0'$ and the general linear group $GL_n$ defined by the formula $\phi (x_{i,j}):=`det(x_{i,j})\neq 0'$. More generally, suppose that $X\subset\bA^n_{\bZ_S}$ is a scheme over $Spec\bZ _S$ given by the equations $f_1(x)=\ldots =f_m(x)=0$, where $f_i(x)\in \bZ _S[x]$. The same equations give us an $\cL _{Rings}(\bZ _S)$-definable set, which we shall also denote by $X$.

Suppose that $U$ and $V$ are $\cT _f$-definable sets given by formulas $\phi (x)$ and $\psi (x)$ respectively, in the same variables. We say that $U$ and $V$ are equal if $\cT _{f}$ contains the sentence $(\forall x)(\phi (x) \leftrightarrow \psi (x))$. It is possible for two non-equal $\cT _f$-definable sets to have the same set of points in some model. However, if two definable sets have the same set of points in every model, then they are equal by the compactness theorem. Similarly, we say that $U$ is contained in $V$ if $\cT _{f}$ contains the sentence $(\forall x)(\phi (x)\to \psi (x))$. The definable sets $V\cap U,V\cup U,V\times U$ are associated with the formulas $\phi (x)\wedge \psi (x) , \phi (x)\vee \psi (x)$ and $\phi (x)\wedge \psi (y)$ respectively, where $y$ is a tuple of variables disjoint from $x$. For the cartesian product, we can omit the requirement that $\phi$ and $\psi$ have the same number of variables.

A $\cT _f$-definable function between the $\cT _f$-definable sets $U$ and $V$ is a $\cT _f$-definable set $W$ that is contained in $U\times V$, such that $\cT _{f}$ implies that $W$ is a graph of a function (note that this can be expressed in $\cL _{Rings}$). A $\cT _f$-definable (linear) group is a $\cT _f$-definable subset $G$ of $GL_n$ such that $\cT _{f}$ implies the axioms of a group for $G$.

Given a $\cT _f$-definable set $X$, the Zariski closure of $X$ is defined in the following way: We look at the ideal of all polynomials $p(x)$ such that $\cT _{f}$ contains the sentence $(\forall x)(\phi (x) \to (p(x)=0))$. This ideal is generated by a finite number of polynomials, say by $p_i(x)$. The Zariski closure of $X$ is the $\cT _f$-definable set given by the formula $p_1(x)=0 \wedge \dots \wedge p_N(x)=0$.

Given a domain $R$, the notions of $\cT _f(R)$-definable sets, functions, and groups are defined similarly. Every $\cT _f$-definable set is a $\cT _f(R)$-definable set, but note that two non-equal $\cT _f$-definable sets can become equal as $\cT _f(R)$-definable sets. For example, the formula $\phi (x):= `1+1=0'$ defines a non empty $\cT _f$-definable set (since it has points over $\bF _2$) but it becomes empty in $\cT _f(\bF_3)$.

We stress again that definable sets are not sets, but rather formulas. The expression ``$x\in V$'' is a synonym for the formula $\phi (x)$ whereas ``$a\in V(L)$'' means that $L$ is a model for our theory, that $a$ is a tuple of elements of $L$, and that $\phi (a)$ holds.

Of course, relative notions are very useful. We will only work over a base which is an affine space, but the definitions can be given for general base variety.

\begin{defn} \begin{enumerate}
\item A $\cT _f$-definable family over $\bA ^n$ is a $\cT _f$-definable subset of $\bA ^n\times \bA ^m$ for some $m$. 
\item A morphism between two definable families $A,B$ over $\bA ^n$ is a definable map $\phi :A \to B$ such that the diagram
\begin{displaymath}
\xymatrix{A \ar[dr] \ar[rr]^{\phi}&  & B \ar[dl] \\ & \bA ^n &}
\end{displaymath}
is commutative.
\item Suppose $A \subset \bA ^n \times \bA ^m$ is a $\cT _f$-definable family defined by the formula $\phi (x,y)$. Given a model $L$ and $a\in \bA ^n (L)$ let $K(a)\subset L$ be the subring generated by the coordinates of $a$. We define the fiber $A_a$ as the $\cT _f(K(a))$-definable set defined by the formula $\phi (a,y)$.
\end{enumerate}
\end{defn}

The fiber product of two definable families over $\bA ^n$ is again a definable family over $\bA ^n$. We denote it by $\times _{\bA ^n}$.

\begin{defn} 
\begin{enumerate}
\item A $\cT _f$-definable family of groups over $\bA ^n$ is a $\cT _f$-definable subset $G$ of $\bA ^n \times GL_m$ such that $\cT _f$ implies that every fiber is a group. 
\item Given a $\cT _f$-definable family of groups $G$ and a $\cT _f$-definable family $\Om$ over the same base, a definable family of actions is a morphism $G \times _{\bA ^n} \Om \to \Om$ such that for every model $L$ and $a\in \bA ^n (L)$, the definable map of the fibers is an action.
\end{enumerate}
\end{defn}

\subsection{Pseudo-finite Fields}
Another theory we shall work with is the asymptotic theory of finite fields, which is also known as the theory of pseudofinite fields of characteristics zero. The language for this theory is again $\cL _{Rings}$. The theory of pseudofinite fields, $\cT _{pf}$, consists of all sentences of $\cL _{Rings}$ that hold for all finite fields, except for the fields of characteristics smaller than $N$ for some $N$. For example, the sentence ``There exists a unique field extension of degree 2, up to isomorphism'' can be expressed in the language of fields. Since it is true for all finite fields, it belongs to $\cT _{pf}$.

\begin{rem} \lbl{rem:use.of.pseudofinite} The use of pseudo-finite fields is for notational simplicity only. If the reader wishes, she can replace all absolute statements of the form ``(the first order sentence) $X$ holds in the theory of pseudo-finite fields'' by the statement ``If $p$ is large enough, then $X$ holds''. 
\end{rem}

Every finite subset of sentences in $\cT _{pf}$ has a model, so by the compactness theorem $\cT _{pf}$ has a model. Note that if $L$ is a model of $\cT _{pf}$ then the characteristics of $L$ is zero (since for every $N$, the theory $\cT _{pf}$ contains the sentences ``The characteristics is different from $N$'').

The notions of $\cT _{pf}$-definable sets, functions etc. are defined similarly. Every $\cT _f$-definable set is a $\cT _{pf}$-definable set. Note, however, that there might be more functions between two definable sets (since the requirement that a set is a graph of a function is stronger in $\cT _f$ than in $\cT _{pf}$).

We denote by $\bN$ the set of non-negative integers. The following theorem is a strengthening of the Lang-Weil estimates (see for example \cite[Theorem 7.1]{Cha-PF} and the references there)

\begin{thm} \lbl{thm:LangWeil} Let $\phi (x,y)$ be a formula in $\cL _{Rings}$. Then there exists a finite set $D\subset \bN \times \bQ _{>0} \cup \{ (0,0)\}$, formulas $\phi _{(d,\mu )} (y)$ for $(d,\mu )\in D$, and a constant $c$, such that the following holds:
\begin{enumerate}
\item The sentence $(\forall y) \bigvee _D \phi _{(d,\mu )}(y)$ holds in the theory of pseudofinite fields.
\item If $p$ is a prime number, $a\in \bF _p ^n$, and $\phi _{(d,\mu )}(a)$ holds, then
\[
\left | |\{ x\in \bF _p ^m | \phi (x,a)\} |-\mu p^{d} \right | < cp^{d-\frac{1}{2}}.
\]
\end{enumerate}
If the Zariski closure of $\phi (x,a)$ is an irreducible variety and has dimension $e$, then $\phi _{(e,\mu )}(a)$ holds for some $\mu$.
\end{thm}

Note that since $(\forall y)\bigvee _D\phi _{(d,\mu )}(y)$ holds in $\cT _{pf}$, then if $p$ is large enough then for every $a\in \bF _p ^n$, there is a $(d,\mu )\in D$ such that $\phi _{(d,\mu )}(a)$ holds.

\begin{defn} A theory $\cT$ is called complete if for any sentence $\phi$, either $\phi \in \cT$ or $\sim \phi \in \cT$. A completion of a theory is a complete theory that contains it.
\end{defn}

By \cite[Theorem 6.14]{Cha-PF}, the completions of $\cT _{pf}$ are given by specifying which integer polynomials are irreducible over the field (and taking all logical implications). This shows that the set of primes for which a single sentence holds is regular in some way.

\begin{defn} \lbl{defn:Artin.set} Let $\cP$ be the set of prime numbers. Given an integer polynomial $f(x)\in \bZ [x]$, let $\cP ^f$ be the set of primes $p$ such that $f(x)$ is irreducible modulo $p$. A set in the Boolean algebra generated by $\cP ^f$ and the Boolean algebra of finite and co-finite sets in $\cP$ is called an Artin set. By the density theorem of Chebotarev, every Artin set is either finite or has a positive analytic density.
\end{defn}

We claim that if $\phi$ is a sentence in the language of fields, then the set of primes $p$ for which $\phi$ holds in $\bF _p$ is an Artin set. For suppose it is not. Enumerate the set of integer polynomials $f_1,f_2,\ldots$ and for each $i$ let $I_i$ be the sentence ``$f_i$ is irreducible''. By our assumption, for every $n$, there are $J_1^n,\ldots ,J_n^n$ such that every $J_i^n$ is either equal to $I_i$ or to $\sim I_i$ and such that both
\[
\cT _{pf} \cup \{ J_1^n \wedge \ldots \wedge J_n^n \wedge \phi \} \quad \textrm{and} \quad \cT _{pf} \cup \{ J_1^n \wedge \ldots \wedge J_n^n \wedge \sim \phi \}
\]
are satisfiable. A diagonalization argument shows that there is a choice $J_1, J_2, \ldots$, where each $J_n$ is equal to either $I_n$ or $\sim I_n$, such that both $\cT _{pf}\cup \{ J_n \} _n \cup \{ \phi \}$ and $\cT _{pf}\cup \{ J_n \} _n \cup \{ \sim \phi \}$ are both satisfyable. But this is a contradiction, since $\cT _{pf} \cup \{ J_n \} _n$ is complete.

\begin{cor} \lbl{cor:size.def} Let $\phi (x,y)$ be a formula in $\cL _{Rings}$. Then there are
\begin{enumerate}
\item A constant $c$.
\item A partition of the set of primes into finitely many Artin sets $\cP _1 ,\dots ,\cP _l$.
\item For each $1\leq i\leq l$, a finite set $D_i \subset \bN \times \bQ _{>0} \cup \{ (0,0)\}$.
\item For each $1\leq i\leq l$, two functions, $e_i:D_i\to\bN$ and $\nu_i:D_i\to\bQ_{>0}$.
\end{enumerate}
such that for every $p\in \cP _i$ and every $a\in \bF _p ^n$, there is a $(d,\mu)\in D_i$ such that
\[
\left | |\{ x\in \bF _q ^m | \phi (x,a)\} |-\mu p^{d} \right | < cp^{d-\frac{1}{2}},
\]
and if we denote by $N_{(d,\mu)}$ the number of the tuples $a\in\bF_p^n$ for which the inequality above holds, then
\[
\left | N_{(d,\mu )}-\nu_i(d,\mu) p^{e_i(d,\mu)} \right | < cp^{e(d,\mu)-\frac{1}{2}}.
\]
\end{cor}

\begin{proof} Let $c,D$ and $\phi _{(d,\mu )}$ be as in Theorem \ref{thm:LangWeil}. For each $(d,\mu )\in D$ apply Theorem \ref{thm:LangWeil} to the formula $\phi _{(d,\mu )}(y)$. In this degenerate case, the theorem says that there are sentences $\phi _{(d,\mu ,e,\nu )}$ such that if $\phi_{(d,\mu ,e,\nu )}$ holds then the number of points in $\phi _{(d,\mu )}(\bF _p)$ is $\nu p^e \pm C\cdot p^{e-\frac{1}{2}}$. Let $\Sigma$ be the set of primes $p$ for which one of the sentences
\[
\bigvee \phi_{d,\mu,e,\nu} \QQ (\forall y)\bigvee \phi_{(d,\mu)}(y)
\]
does not hold. Since these sentences hold in $\cT_{pf}$, we get that $\Sigma$ is finite. By the above, there is a partition of the primes into Artin sets $\cP_i$ such that for each $i$ and $(d,\mu,e,\nu)$, the sentence $\phi_{(d,\mu,e,\nu)}$ holds for all $\{\bF_p | p\in \cP_i\}$ or for none. By further partitioning of the $\cP_i$, we can assume that for each $i$ either $\cP_i$ is infinite and $\cP_i\cap\Sigma=\emptyset$, or $\cP_i$ is a singleton. For each $i$ such that $\cP_i\cap\Sigma=\emptyset$, set $D_i=\{(\mu,d) | (\exists y) \phi_{(d,\mu)}(y)$ and let $(e(d,\mu),\nu(d,\mu))$ be the unique $(e,\nu)$ such that $\phi_{(d,\mu,e,\nu)}$ holds for all the primes in $\cP_i$. For $\cP_i$ a singleton, set $D_i=\{(0,1)\}, e(0,1)=\nu(0,1)=1$. It is clear that if $c$ is large, then the proposition holds.
\end{proof}

\subsection{Definable Families of Groups} 

\begin{prop} \lbl{prop:def.radical} Let $L\subset\bA^n\times\M_n$ be a $\cT_f$-definable family of Lie algebras. Then there is a definable family $R\subset L$ such that for every model $F$ of $\cT_f$ of high enough characteristics, and for every $x\in\bA^n(F)$, the fiber $L_x$ is the unipotent radical of the Lie algebra $L_x$.
\end{prop}

\begin{proof} It is known that given $n$, if $F$ is a field of high enough characteristics, and $\sL\subset\M_n(F)$ is a Lie algebra, then an element $x\in\sL$ is in the unipotent radical of $\sL$ if and only if for every $y\in\sL$, the element $[x,y]$ is nilpotent.
\end{proof}

For a root datum $\Phi$ denote $H _{\Phi}$ the adjoint algebraic group attached to $\Phi$.

\begin{lem} \lbl{lem:def.root.datum} Let $G\subset \bA ^m \times GL_n$ be a $\cT _{pf}$-definable family of semisimple adjoint groups. Then there is a definable partition $\bA ^m = X _1 \sqcup \ldots \sqcup X _k$, and for every $i$ there are root data $\Phi _{i,1},\ldots ,\Phi _{i,j}$, such that for $p$ large enough, and for $x\in X_i(\bF_p)$, the group $G_x(\bF_p)$ is isomorphic to $H _{\Phi _{i,1}}(\bF_p) \times \ldots \times \uH _{\Phi _{i,j}}(\bF_p)$.
\end{lem}

\begin{proof} We first assume that $G$ is a family  of simple adjoint groups. In this case we need to show that for every root datum $\Phi$, the set of $x\in \bA ^m$ such that $G_x$ is isomorphic to $\uH _{\Phi}$, is a definable set. A-priori, this condition is not definable, as the isomorphism can be a polynomial map of very high degree.

Let $\sL _{\Phi}$ be the Lie algebra attached to $\Phi$. Define a family of Lie algebras $L \subset \bA ^m \times {gl_n}$ as follows: for $x\in \bA ^m$, let $U _x \subset G _x$ be the set of unipotent elements. Define $L _x$ to be the span of $\log U_x$.

It is known that $G_x(\bF_p)$ is isomorphic to $\uH _{\Phi}(\bF_p)$ if and only if $L _x(\bF_p)$ is isomorphic to $\sL _{\Phi}$. This, however, is a definable condition, since every morphism between Lie algebras is linear.

The argument for products of simple groups is similar.
\end{proof}

\begin{rem} This proof actually shows that every map between connected algebraic groups $G \subset {GL_n}$ and $H\subset {GL_m}$ can be represented by a polynomial whose degree is bounded as a function of $m$ and $n$.
\end{rem}

\begin{prop} \lbl{prop:filtration.subgroups} Let $X$ be a definable set in $\cT_{pf}$, and let $S \subset X\times {GL_n}$ be a family of definable groups in $\cT_{pf}$. Then there are:
\begin{enumerate}
\item A definable partition $X=X_1\sqcup \ldots \sqcup X_m$,
\item For each $i$, a finite sequence of root data $\Phi^i_1,\ldots,\Phi^i_{n_i}$,
\item Definable families $U \subset C \subset S$ of normal subgroups, 
\item A constant $c$,
\end{enumerate}
such that for every $p$ large enough, for any $i$, and for any $a\in X_i(\bF _p)$, the following hold:
\begin{enumerate}
\item \lbl{cond:U} $U_a(\bF_p)$ is a unipotent group. Moreover, $U _a(\bF _p)$ is the maximal normal $p$-subgroup of $S _a(\bF_p)$.
\item \lbl{cond:C.by.U} $C _a (\bF _p)/U _a(\bF _p)$ is isomorphic to $H_{\Phi^i_1}(\bF _p) \times \ldots \times H_{\Phi^i_{n_i}}(\bF _p)$.
\item \lbl{cond:S.by.C} The group $S _a(\bF _p)/C _x(\bF _p)$ is an extension of an abelian group, whose order is prime to $p$, by a group of size less than $c$.
\end{enumerate}
\end{prop}

\begin{proof}  Let $P \subset X\times {GL_n}$ be the definable family
\[
P = \{ (x,g)\in S | (g-1)^n=0 \} .
\]
Since on unipotent elements, the function $\log$ is a polynomial, the family
\[
L=\{ (x,A)\in X\times {\mathfrak gl}_n | \textrm{$A$ is in the span of $\log P_x$} \}
\]
is a definable family. It is easy to see (see \cite[Lemma 1.6]{No}) that if $p$ is large enough, then for every $x\in X(\bF _p)$, the set $L_x(\bF _p)$ is a Lie algebra. By Proposition \ref{prop:bounded.generation.algebraic}, there is a family $C\subset X\times GL_n$ such that for every $p$ large enough, and every $x\in X(\bF_p)$, the set $C_x(\bF_p)$ is equal to the group generated by the set $\exp L_x(\bF_p)$. By \cite[Theorem B]{No}, we have that $C_x(\bF_p)\subset S_x(\bF_p)$ for all $p$ large enough, and hence $C\subset S$. Clearly, $C$ is a family of normal subgroups of S. By \cite[Theorem C]{No}, for every $p$ large enough and every $x\in X(\bF _p)$, there is a commutative subgroup $H\subset S_x(\bF _p)$ such that $HC_x(\bF _p)$ is normal in $S_x(\bF _p)$ and its index is less than a constant, which we denote by $c$. Moreover, by the same theorem, the order of $H$ is prime to $p$. Claim (\ref{cond:S.by.C}) follows from this.

By Proposition \ref{prop:def.radical}, there is a definable family of Lie subalgebras $L^u \subset L$, such that for every $x\in X$, the Lie subalgebra $L^u _x$ is the unipotent radical of $L_x$. By Proposition \ref{prop:bounded.generation.algebraic}, there is a definable family $U\subset C$ such that for every $p$ large enough, and every $x\in X(K)$, the set $U_x(\bF_p)$ is the subgroup generated by the set $\exp L^u_x(\bF_p)$. Clearly, $U_x(\bF_p)$ is a unipotent normal $p$-subgroup of $C_x(\bF_p)$. Since $U_x(\bF _p)$ is characteristic in $C_x(\bF _p)$ and since, if $p$ is large enough, the order of $S_x(\bF _p)/C_x(\bF _p)$ is not divisible by $p$, we get that $U_x(\bF _p)$ is a normal $p$-subgroup of $S_x(\bF _p)$. This finishes the proof of claim (\ref{cond:U}), except for the maximality.

For every $p$ large enough and for every $x\in X(\bF_p)$, the Lie algebra $L_x(\bF_p)/L_x^u(\bF_p)$ is reductive. Since it is generated by nilpotents, it is, in fact, semisimple. Theorem $B$ of \cite{No} shows that $C_x(\bF_p)/U_x(\bF_p)$ is equal to $H(\bF_p)^+=[H(\bF_p),H(\bF_p)]$, where $H$ is a semisimple algebraic group with Lie algebra $L_x(\bF_p)/L_x^u(\bF_p)$. By lemma \ref{lem:def.root.datum} we get a definable partition $X=X_i \sqcup \ldots \sqcup X_m$ and root data $\Phi^i _j$ such that for every $p$ large enough and every $x\in X_i(\bF _p)$, the group $C_x (\bF_p)/U_x(\bF_p)$ is isomorphic to
\[
H_{\Phi ^i _1}(\bF_p) \times \ldots \times H_{\Phi ^i _{n_i}}(\bF_p).
\]
Claim (\ref{cond:C.by.U}) follows. Also, it follows that $C_x (\bF_p)/U_x(\bF_p)$ does not have any normal $p$ groups, and hence the maximality claim in (\ref{cond:U}) follows.
\end{proof}

\subsection{Henselian Valued Fields}
We will also work in the theory of valued Henselian fields. The language for this theory is the language $\cL _{Vf}$ of valued fields, which we proceed to define. In $\cL _{Vf}$ one can quantify over three kinds of variables (which are called sorts): one is the valued fields sort, one is the residue field sort, and one is the value group sort. The language has constants $0,1$ (which are valued field sort), $\bar{0},\bar{1}$ (which are residue field sort) and $\tilde{0}$ (which is value group sort). The relations are equality and $\tilde{<}$, but we can only equate expressions of the same sort (so $0=\bar{0}$ is not a legitimate expression) and only $\tilde{<}$-compare expressions which are in the value group sort. The functions symbols are two sets of addition and multiplication (for the value field and residue field sorts), addition (for the value group sort), and two additional function symbols: a function $\val$ (called valuation) from the valued field sort to the value group, and a function $\ac$ (called angular component) from the value field sort to the residue field sort. Here again, there are definable sets, functions etc. An example of a definable set is $\cO$, which is defined by the formula $\phi(x):=`\val(x)\geq 0'$, where $x$ is a variable of valued field sort. Since in $\cL _{Vf}$ there is more than one sort, there might be confusion regarding the variables of the definable sets. We resolve this confusion by adding the subscripts $V,R,O$ for valued field, residue field and value group respectively. So, for example $\bA ^n _V$ is the affine space whose coordinates are valued field sort and $(GL_n)_R$ is the set of invertible $n\times n$ matrices whose entries are from the residue field.

The theory $\cT _{Hvf}$ of Henselian valued fields, whose valuation group is elementary equivalent to $\bZ$, consists of the following axioms:
\begin{itemize}
\item The axioms of fields for the valued field sort and for the residue field sort.
\item The axioms of non archimedian valuation.
\item All sentences that hold for $\bZ$ for the value group sort.
\item The sentences $(\forall x,y\neq0)\ac(xy)=\ac(x)\cdot\ac(y)$, $(\forall x,y\neq0) (\val(x)<\val(y) \to \ac(x+y)=\ac(x)$, and $(\forall x,y\neq 0)(\val(x)=\val(y) \wedge \ac(x)\neq -\ac(y))\to \ac(x+y)=\ac(x)+\ac(y)$.
\item Sentences stating that the field is Henselian. 
\end{itemize}

We shall be mainly interested in the models $\bM _p$ of $\cT _{Hvf}$, which interpret $\bA_V$ as $\bQ _p$, interpret $\bA_R$ as $\bF _p$, interpret $\bA_O$ as $\bZ$, interpret $\val(x)$ as the $p$-adic valuation of $x$, and interpret $\ac(x)$ as the first non zero coefficient in the $p$-adic expansion of $x$.

Let $\cT_{Hvf,0}$ be the theory $\cT_{Hvf}$ together with the axioms that claim that the characteristic of the residue field is equal to zero. While no $\bM_p$ is a model for $\cT_{Hvf,0}$, it follows from the compactness theorem that every sentence that holds in $\cT_{Hvf,0}$ is also true in all but finitely many of the $\bM_p$'s. 

\begin{thm} (Elimination of quantifiers in Henselian fields---\cite[Theorem 4.1]{Pa}) \lbl{thm:elim.quan.valued} Let $\phi (x,y,z)$ be a formula in the language $\cL _{Vf}$ where the variable $x$ is of the valued field sort, the variable $y$ is of residue field sort, and the variable $z$ is of value group sort. Then there is a partition of $\bA _F ^n$ into finitely many constructible sets $C_j$; and for each $j$ there are polynomials $P_1 ^j(x),\dots ,P_n ^j(x)$ that do not vanish on $C_j$, formulas $\psi _1 ^j(x,y)$ in the language of rings, and formulas $\psi _2 ^j(x,y)$ in the language of ordered groups, such that $\cT _{Hvf,0}$ implies that $\phi (x,y,z)$ is equivalent to the formula
\[
\bigvee _j \left ( x\in C_j \wedge \psi _1 (\ac(P_i (x)),y) \wedge \phi _2 (\val(P_i(x)),z) \right )
\]
\end{thm}

\begin{thm} (Elimination of quantifiers in the theory of $\bZ$---\cite[Lemma 5.5]{Pa}) \lbl{thm:elim.quan.ordered} Let $A$ be a definable set in the language of ordered groups. Then the theory of $\bZ$ implies that $A$ is equal to a Boolean combination of definable sets defined by formulas of the form
\[
\phi (x) \geq 0 \wedge (\exists y\in \Ga) \psi (x)=n\cdot y
\]
where $\phi (x),\psi (x)$ are linear functionals with integer coefficients and $n\in \bN$.
\end{thm}

The following is an easy corollary of Theorem \ref{thm:elim.quan.ordered}:
\begin{lem} \lbl{lem:ord.def.func} Every definable function $f:\bA_O^n\to\bA_O$ is piecewise linear. In other words, there is a partition of $\bA_O^n$ to definable sets $\bA_O^n=X_1\sqcup\ldots\sqcup X_m$, there are linear functionals $\varphi_1,\ldots,\varphi_m$ with rational coefficients, and there are definable elements $\ga_1,\ldots,\ga_m\in\bA_O$ such that $\cT_{Hvf}$ implies the sentence
\[
(\forall x\in\bA_O^n) \left ( x\in X_i \to f(x)=\varphi_i(x)+\ga_i \right ).
\]
\end{lem}

\begin{prop} \lbl{prop:func.to.ordered} A definable function $q:\bA_V ^n\to \bA_O^1$ is (in $\cT_{Hvf,0}$, and, a-posteriori, in $\bM_p$ for $p$ large enough) piecewise of the form $\frac{1}{n} val\left ( \frac{f(x)}{g(x)} \right )$ where $f$ and $g$ are polynomials.
\end{prop}

\begin{proof} Let $q:\bA^n_V \to \bA_O^1$ be a definable function. By Theorem \ref{thm:elim.quan.valued}, the graph of $q$, which is a subset of $\bA_V^n\times\bA_O^1$, is defined by a formula of the type
\[
\bigvee \phi _i (\val(P_1(x)),\dots ,\val(P_r(x)),\ga ) \wedge \psi _i (\ac(P_1(x)),\dots ,\ac(P_r(x))) \wedge (Q^i_1(x) \dots = Q^i_t(x)=0)\wedge(Q_{t+1}^i(x)\neq0)
\]
where $P_i(x)$ and $Q^i _j(x)$ are polynomials, $\phi (y_1,\dots ,y_r,z)$ is a formula in the language of ordered groups, and $\psi (x_1,\dots ,x_r)$ is a formula in the language of fields. Decompose the domain of $q$ according to the conditions $Q^i _j (x)=0$ and $\psi _i(\ac(P_j(x)))$. Let $A$ be one of the pieces. On $A$, the graph of $q$ is given by a formula $\phi (\val(P_j(x)),\ga)$. Again, by Theorem \ref{thm:elim.quan.valued}, $\val(A)$ is a subset $B$ in $\bA_O^n$ which is definable in the language of ordered groups. Therefore the formula $\phi (y_1,\dots ,y_r,z)$ defines a graph of a function from $B$ to $\bA_O^1$. Since by Lemma \ref{lem:ord.def.func}, every such function is piecewise linear, we get that after a further division of the domain, $q$ is of the form required.
\end{proof}

\begin{prop} \lbl{prop:red.def.fam} Let $A \subset \bA _V ^n \times \bA _R ^m$ be a definable family. Then there is a definable set $B\subset \bA _R ^{l+m}$ and a definable function $f:\bA _V ^n\to \bA _R ^l$ such that 
\[
A=\{ (x,y)\in\bA_V^n\times\bA_R^m\qq|\qq(f(x),y)\in B\}
\]
in $\cT_{Hvf,0}$, and, a-posteriori, in all but finitely many of the models $\bM_p$. We say that $A$ is the pullback of $B$ via $f$.
\end{prop}

\begin{proof} By elimination of quantifiers we can assume that $A$ is defined by a formula $\phi (x,y)$ ($x$ is a valued field sort and $y$ is a residue field sort) that is the conjunction of conditions of the form 
\begin{enumerate}
\item $\psi (\val(P_1(x)),\dots ,\val(P_n(x)))$, where $P_i$ are polynomials and $\psi$ is a formula in the language of ordered group.
\item $\xi (\ac(Q_1(x)),\dots ,\ac(Q_m(x)),y)$, where $Q_i$ are polynomials and $\xi$ is a formula in the language of fields.
\end{enumerate}
Decompose $\bA_V^n$ according to condition (1) Denote the resulting pieces by $X_i$, where $i\in I$. For every $i$ there is a formula $\xi _i$ such that the restriction of $A$ over $X_i$ is the pullback of the definable set $\xi _i$ via the map $(\ac(q_1),\dots ,\ac(q_m))$. Define a map $\Psi :\bA_V^n\to\bA_R^{m+|I|}$ by $\Psi (x)_i=\ac(q_i(x))$ for $i\leq m$, $\Psi (x)_j =1$ if $j>m$ and $x\in X_{j-m}$, and $\Psi (x)_j =0$ if $j>m$ and $x \not \in X_{j-m}$. Let $\Om \subset \bA_R^{m+|I|}$ be the definable set that is the conjugation of the conditions
\[
z_{j+m}=1 \to \xi _j (z_1 ,\dots ,z_m)
\]
for $j=1,\dots ,|I|$. It is now clear that $A$ is the pullback of $\Om$ via the map $\Psi$.
\end{proof}

\begin{prop} \lbl{prop:def.polynom} \begin{enumerate}
\item Given a polynomial $P(x)\in \bQ_{\geq0}[x]$ such that $\lim_{x\to\infty}P(x)=\infty$, there is a definable set $Y$ in the language of rings and a constant $c$ such that for all primes $p$,
\[
1-c\cdot p^{-\frac{1}{2}} < \frac{|Y (\bF _p)|}{P(p)} <1+c\cdot p^{-\frac{1}{2}} .
\]
\item Let $X\subset\bA_V^n$ be a definable set, and let $A$ be a definable family over $X$. Then there is a definable function $\psi :X\to \bA_O$, a definable family $B$ over $X$, and a constant $c$, such that for all primes $p$ and $x\in X(\bM _p)$, either $A _x (\bM _p)$ is empty or
\[
1-c\cdot p^{-\frac{1}{2}} < p^{\psi (x)} \cdot |B _x (\bM _p)| \cdot |A _x(\bM _p)| <1+c\cdot p^{-\frac{1}{2}} .
\]
\end{enumerate}
\end{prop}

\begin{proof} (1) Suppose that the leading coefficient of $P(x)$ is $\frac{a}{b}x^k$ where $a,b,k$ are positive integers. Let $\phi :C_1 \to C_2$ be a Galois cover defined over $\bQ$ of irreducible curves with Galois group $\bZ /b$. For almost all primes $p$, the reduction modulo $p$ of $\phi$ is also a Galois cover with the same Galois group. Let $D$ be the definable set defined by the formula
\[
x\in C_2 \wedge (\exists y\in C_1)(x=\phi (y)) .
\]
Then by Weil's theorem there is a constant $K$ such that 
\[
1-K\cdot p^{-\frac{1}{2}} < \left | \frac{|D(\bF _p)|}{\frac{1}{b}p} \right | <1+K\cdot p^{-\frac{1}{2}} .
\]
To get the claim of the lemma, take $\cP$ be the definable set
\[
(\underbrace{D\sqcup \dots \sqcup D} _{a} )\times \bA ^{k-1}.
\]

(2) By Proposition \ref{prop:red.def.fam} and Theorem \ref{thm:LangWeil} there is a constant $K$ and a partition of $X$ into definable sets $X_i$, and for each $X_i$ there is $d_i\in \bN$ and $c_i\in \bQ _{>0}$, such that if $x\in X_i (\bM _p)$, then
\[
\left | |\cA _x (\bM _p)|-c_ip^{d_i} \right | < K\cdot p^{d_i-\frac{1}{2}}.
\]
Using the construction from (1), one can find definable sets $B_i$ such that for all $p$,
\[
1-K\cdot p^{-\frac{1}{2}} < \left | \frac{|B_i(\bF _p)|}{c_i \cdot p} \right | < 1+K\cdot p^{-\frac{1}{2}} .
\]
Denote by $\cB$ the definable family that is equal to $B_i$ over $X_i$, and denote by $\psi _i:X\to \bA_O^1$ the definable function that equals $d_i-1$ on $X_i$. Then $\cB ,\psi$ satisfy the requirements of the lemma.
\end{proof}

\subsection{$V$-Functions}

Definable elements in $\bA_O^1$ give us a collection of (rational) numbers, indexed by the prime numbers: Given a definable element $\ga\in\bA_O^1$ and a prime $p$, we consider the number $\ga_p=p^{\ga^{\bM_p}}$, where $\ga^{\bM_p}\in\bZ$ is the interpretation of $\ga$ in the model $\bM_p$. More generally, definable functions from $\bA_V^n$ to $\bA_O^1$ give us a collection of rational-valued functions. Namely, if $f:\bA_V^n\to\bA_O^1$ is a definable function, and $p$ is a prime number, we consider the function $f_p:\bQ_p^n\to\bQ$ given by $x\mapsto p^{f^{\bM_p}(x)}$, where $f^{\bM_p}$ is the interpretation of $f$ in the model $\bM_p$.

Another source of numbers in $\cT_{Hvf}$ is definable sets in $\bA_R^m$. Given a definable set $X\subset\bA_R^m$ and a prime $p$, we consider the number $X_p=|X(\bM_p)|$ (since $X(\bM_p)\subset \bA_R^m(\bM_p)=\bF_p^m$, the set $X(\bM_p)$ is indeed finite). As above, from a definable set $Y\subset\bA_V^n\times\bA_R^m$ we get a collection of integer-valued functions, indexed by the prime numbers. The next definition is a generalization of these two constructions:

\begin{defn} \lbl{defn:V-function}
Let $X\subset \bA _V ^n$ be a definable set.
\begin{enumerate}
\item A $V$-function with domain $X$ is a tuple of the form $\cF=(X_i,\phi _i ,\cV _i)_{i\in I}$, where $I$ is a finite set and
\begin{enumerate}
\item $X_i\subset \bA_V^n$ are definable sets that form a partition of $X$.
\item $\phi _i:\bA_V^n \to \bA^1_O$ are definable maps, and 
\item $\cV _i\subset \bA_V^n \times \bA_R^{n_i}$ are definable sets.
\end{enumerate}

\item Given a $V$-function $\cF$ with domain $X$ and a prime number $p$, we define a function $\cF _p :X(\bM _p)\times \bC \to \bC$ by
\[
\cF _p (x,s)=\sum _{i\in I} 1_{X_i(\bM _p)}(x) p^{-s\phi _i (x)}\cdot | \cV _i (\bM _p)_x|^{-s} .
\]

\item A $V$-function with domain $X$ is called bounded if there exists a definable element $\ga\in\bA_O$ such that for all $i$, the following sentences hold
\begin{enumerate}
\item $\forall x\in X_i\qq(\phi_i(x)>\ga)$.
\item $\forall(x_1,\ldots,x_n)\in X_i\left ( (\val(x_1)<\ga\wedge\ldots\wedge\val(x_n)<\ga)\longrightarrow(\cV_i)_x=\emptyset\right )$.
\end{enumerate}
\end{enumerate}
\end{defn}

\begin{rem} 
\begin{enumerate}
\item By changing $X_i$, every sequence of functions of the more general form
\[
(x,s)\mapsto\sum _{i\in I} 1_{X_i (\bM _p)}(x) p^{-s\phi _i (x)+\psi _i(x)}\cdot | \cV _i (\bM _p)_x|^{-s} \cdot |\cW _i (\bM _p)_x| , 
\]
where $X_i\subset\bA_V^n$ are definable sets, $\phi _i,\psi _i : X_i \to \bA_O^1$ are definable maps, and $\cV_i,\cW _i \subset \bA _V ^n\times\bA _R^m$ are definable families, actually comes from a $V$-function.

\item If $\cF _1$ and $\cF _2$ are $V$-functions, then there are $V$-functions $\cG _{+}$ and $\cG _{\times}$ such that for every $p$
\[
(\cG _{+})_p=(\cF_1)_p +(\cF _2)_p
\]
and
\[
(\cG _{\times})_p=(\cF_1)_p \cdot(\cF _2)_p.
\]
We shall write $\cF _1 +\cF _2$, respectively $\cF _1 \cdot \cF_2$ instead of $\cG _+$, respectively $\cG _\times$.
\end{enumerate}
\end{rem}

\begin{exam} \lbl{exam:geometric.series} Fix integers $A$ and $B$. Let $\cO^\times\subset\bA^1_V$ be the definable set given by the formula `$x\neq0\wedge\val(x)\geq0$', let $\phi:\cO^\times\to\bA_O^1$ be the function $\phi(x)=A\val(x)$, and $\psi:\cO^\times\to\bA_O^1$ be the function $\psi(x)=(B+1)\val(x)$. By the previous remark we get a $V$-function $\cF$ such that
\[
\cF_p(x,s)=\left\{ \begin{matrix} p^{n(As+B+1)} & \val(x)=n\geq0\\0 & \textrm{else} \end{matrix}\right.
\]
and so
\[
\int\cF_p(x,s)dx=\frac{p-1}{p}\sum_{n=0}^\infty p^{As+B}=\frac{p-1}{p}\frac{1}{1-p^{As+B}}
\]
\end{exam}

\section{Uniformity of the Local Factors I} \lbl{sec:Uniformity.I}

Our goal in this section is to prove

\begin{thm} \lbl{thm:uniformity} Let $\Sigma$ be a finite set of primes, and let $\uG$ be a linear algebraic group scheme over $Spec\bZ _\Sigma$ such that the generic fiber of $\uG$ is semisimple, simply connected, and connected. There is a definable set $X\subset \bA ^n _F$ and a $V$-function $\cF$ with domain $X$, such that the sequence of functions $\ze _{\uG (\bZ _p )}(s)-1$ is equivalent to the sequence of functions
\[
\xi _p (s) = \int _{X(\bM _p)}\cF _p (x,s) d\la (x)
\]
where $d\la$ is the restriction of the Haar measure of $\bQ _p ^n$
to $X(\bM _p)$.
\end{thm}

By Theorem \ref{thm:jaikin}, Proposition \ref{prop:def.polynom}, and Example \ref{exam:geometric.series}, it is enough to show that there is a prime $p_0$ and a $V$-function $\cF$ such that the sequence $(\ze_{\uG(\bZ_p)}(s)-1)_{p>p_0}$ is equivalent to the sequence $(\xi_p(s))_{p>p_0}$ considered in the theorem. Therefore, in this section we will assume that $p$ is large enough.

\subsection{Representations of the First Congruence Subgroup} \lbl{subsec:reps.first.congruence}

Let $\uG$ be as in theorem \ref{thm:uniformity}. The corresponding $\cT_{Hvf}$-definable group will be denoted by $G$. The definable subset $G_\cO$ is defined by
\[
g=(g_{i,j})\in G_\cO\iff g\in G\wedge\val(g_{i,j})\geq 0.
\]
For all $p\not\in\Sigma$ we have that $G_\cO(\bM_p)=\uG(\bZ_p)$, a group which was denoted by $G_p$.

Let $\fg \subset (\M_n)_F$ be the definable set such that, for almost all $p$'s, $\fg (\bM_p)$ is the Lie algebra $\fg _p$ of $G_p$; see Subsection \ref{subsec:Lie.algebras} for the construction. The definable set
\[
\fg ^1 = \{ A\in \fg | \val(A_{i,j})>0 \}
\]
satisfies that $\fg ^1(\bM_p)$ is the Lie algebra $\fg_p^1$ of $G_p^1$---the first congruence subgroup of $G_p$---for almost all $p$'s.

Assume that $p$ is large. The orbit method (Theorem \ref{thm:orbit.method}) gives us a map $\Xi _p$ from $(\fg _p^1) ^{\vee}$ onto $\Irr(G_p ^1)$ such that $\Xi _p (\theta )=\Xi _p (\theta ')$ if and only if there is a $g\in G_p ^1$ such that $Ad^*(g)\theta =\theta '$. 

Recall that $\cO\subset \bA _V ^1$ is the definable set attached to the formula `$\val(x)\geq 0$'. Let $\X$ be the following definable set
\[
\X = \{ A\in \fg _F\qq|\qq \max \{ \val(A_{i,j})\} =0 \} \times \cO .
\]
Consider the function $x\mapsto \exp (2\pi i x)$ defined on $\bZ[\frac{1}{p}]$. It has a unique extension to a continuous character of $\bQ _p$, which we also denote by $\exp (2\pi ix)$. Let $<,>$ be the bilinear form 
\[
<A,B>=\trace(A\cdot B)
\]
on the space of matrices.

For every prime $p$, the map $\Phi _p :\X (\bM _p)\to (\fg _p ^1)^{\vee}$ given by
\[
\Phi _p ((A,z))(B)=\exp \left ( \frac{2\pi i}{z} <A,B> \right )
\]
is a surjection. Let $\Psi _p$ be the composition of $\Xi _p$ and $\Phi _p$. We have the following

\begin{thm} \lbl{thm:jaikin.precise} (\cite[Theorem 4.6]{Jai}) There are definable functions $\phi _1,\phi _2 :\X \to \bA_O ^1$ such that for all primes $p$, and for all $(A,z)\in \X(\bM _p)$,
\[
\dim\Psi _p (A,z)=p^{\phi _1(A,z)}
\]
and
\[
\la(\Psi _p ^{-1}(\Psi _p (A,z)))=p^{\phi _2 (A,z)}.
\]
\end{thm}

\subsection{Decomposition Trees} We describe our method for extending representations from the first congruence subgroup, $G_p^1$, to the whole group, $G_p$. Let $\rho$ be an irreducible representation of $G_p^1$. Recall that $\Irr(G_p|\rho)$ is the set of irreducible representations of $G_p^1$, whose restrictions to $G_p^1$ contain $\rho$ as a component. We wish to compute the relative zeta function
\[
\ze_{G_p|\rho}(s) = \sum _{\tau \in \Irr(G_p|\rho)} \left ( \frac{\dim \tau}{\dim \rho} \right ) ^{-s} .
\]

Consider the stabilizer $G_p^1 \subset S \subset G_p$ of the representation $\rho$. Let $V$ be the maximal normal pro-$p$ subgroup of $S$.

\begin{lem} \lbl{lem:recursive.rel.zeta} Let $K\subset V \subset S \subset H$ be an increasing chain of groups, and let $\rho$ be an irreducible representation of $K$. Assume that $K$ is normal in $H$, that $V$ is normal in $S$, that $S$ is the stabilizer of $\rho$ in $H$, and that $K$ is of finite index in $H$. For a representation $\tau$ of $V$, denote the orbit of $\tau$ under the action of $S$ by $\tau^S$. Then
\[
\ze_{H|\rho}(s) = \sum _{\tau \in \Irr(V|\rho)} \frac{[H:S]^{-s}}{|\tau ^S|} \left ( \frac{\dim \tau}{\dim \rho} \right ) ^{-s} \ze_{S|\tau}(s).
\]
\end{lem}

\begin{proof} Let $\chi \in \Irr (S|\rho)$. The irreducible components of the restriction $\Res^S_K\chi$ form one $S$-orbit, since $\chi$ is irreducible and $K$ is normal in $S$. Since this restriction contains $\rho$ as an irreducible component, we deduce that $\Res^S_K\chi =\Res^V_K \circ \Res^S_V \chi$ is a multiple of $\rho$. Therefore, every irreducible component of $\Res^S_V\chi$ is in $\Irr(V|\rho)$. 

Let $\tau \in \Irr(V|\rho)$. The irreducible components of $\Res ^S_V\Ind^S_V \tau$ are just $\tau ^S$. Therefore, if $\tau_1,\tau_2$ are two representations in $\Irr(V|\rho)$, which are in the same $S$-orbit, then $\Irr(S|\tau_1)=\Irr(S|\tau_2)$, whereas if $\tau_1$ and $\tau_2$ are not in the same orbit, then $\Irr(S|\tau_1)$ and $\Irr(V|\tau_2)$ are disjoint.

Let $\tau_1,\ldots,\tau_m$ be representatives for the $S$-orbits in $\Irr(V|\rho)$. We compute:
\[
\ze_{S|\rho}(s) = \sum _{\chi \in \Irr(S|\rho)} \left ( \frac{\dim \chi}{\dim \rho} \right ) ^{-s}= \sum _{i=1}^m \sum _{\chi \in \Irr(S|\tau_i)} \left ( \frac{\dim \chi}{\dim \rho} \right ) ^{-s} = 
\]
\[
= \sum _{\tau \in \Irr(V|\rho)} \frac{1}{|\tau ^S|} \left ( \frac{\dim \tau}{\dim \rho} \right ) ^{-s} \sum _{\chi \in \Irr(S|\tau)} \left ( \frac{\dim \chi}{\dim \tau} \right ) ^{-s} = \sum _{\tau \in \Irr(V|\rho)} \frac{1}{|\tau ^S|} \left ( \frac{\dim \tau}{\dim \rho} \right ) ^{-s} \ze_{S|\tau}(s).
\]
Since for every representation $\chi \in \Irr(S|\rho)$, the induction $\Ind^H_S\chi$ is irreducible, and since all the irreducible representations of $H$, lying over $\rho$, are obtained in this way, we get that
\[
\ze_{H|\rho}(s) = \ze_{S|\rho}(s) \cdot [H:S]^{-s}.
\]
\end{proof}

Lemma \ref{lem:recursive.rel.zeta} reduces the computation of $\ze_{G_p|\rho}(s)$ to the computation of $\Irr(V|\rho)$, and, for each $\tau \in \Irr(V|\rho)$, a computation of $\ze_{S|\tau}(s)$.

Let $\Irr(V|\rho)=\{ \rho_1,\ldots ,\rho _n\}$. For each $\rho_i$, let $S_i$ be the stabilizer of $\rho_i$ in $S$, let $V_i$ be the maximal normal pro-$p$ subgroup of $S_i$, and let $\rho_{i,1},\ldots ,\rho_{i,n_i}$ be the irreducible characters of $V_i$ lying over $\rho _i$. The following diagram is a summary of the notation so far:

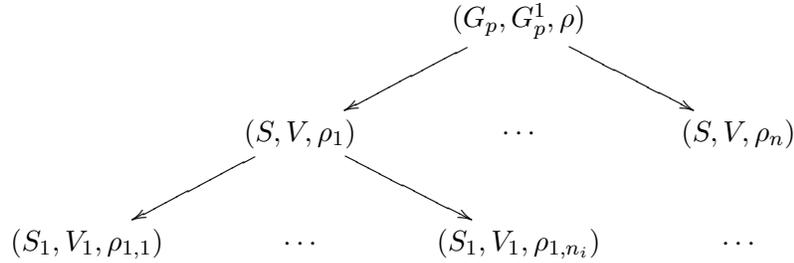
\begin{displaymath}
\xymatrix{& & (G_p,G_p^1,\rho) \ar[dr] \ar[dl] & \\ & (S,V,\rho_1) \ar[dr] \ar[dl] & \cdots & (S,V,\rho_n) \\ (S_1,V_1,\rho_{1,1}) & \cdots & (S_1,V_1,\rho_{1,n_i}) & \cdots}
\end{displaymath}
We may continue in the same fashion, constructing a deeper and deeper trees. We reach a leaf of the tree whenever $S_{i_1\cdots_k}$ is the stabilizer of $\rho_{i_1\cdots i_{k+1}}$ and $S_{i_1\cdots i_k}/V_{i_1\cdots i_k}$ has no non-trivial normal $p$-subgroups.

We call the resulting tree, whose vertices are labeled by triples $(S_{i_1\cdots i_k},V_{i_1\cdots i_k},\rho_{i_1\cdots i_{k+1}})$, the {\em decomposition tree} of $\rho$. The relative representation zeta function $\ze_{G_p|\rho}(s)$ can be easily computed from the decomposition tree, and the relative representation zeta functions of the leaves. The following lemma shows that the zeta functions of the leaves are simple:

\begin{lem} \lbl{lem:zeta.leaves} For every $n$, there is a constant $c$, that depends only on $n$, such that if $p$ is a prime number, which is large enough with respect to $n$, then the following is true: Let $S\subset \GL_n(\bZ_p)$ be a group, let $V$ be a normal pro-$p$ subgroup of $S$ that contains the first congruence subgroup $G_p^1$, and let $\rho$ be a representation of $V$. Assume that $S$ stabilizes $\rho$, and assume that $S/V$ has no non-trivial normal $p$-subgroups. Then
\begin{enumerate}
\item The group $(S/V)^+$ is a perfect extension of a direct product of finite simple groups of Lie type.
\item $c^{-1-s} \cdot \ze_{S|\rho}(s) \leq [S:S^+] \cdot \ze_{(S/V)^+}(s) \leq c^{1+s} \cdot \ze_{S|\rho}(s).$
\end{enumerate}
In the middle term, $\ze_{(S/V)^+}(s)$ is the (non-relative) representation zeta function of the group $(S/V)^+$, and $S^+$ is the closed subgroup of $S$ that is generated by the pro-$p$ elements of $S$.
\end{lem}

\begin{proof} Denote the quotient $S/G_p^1$ by $\Ga$. By a theorem of Larsen and Pink (theorem 0.2 of \cite{LP}), there are normal subgroups $\Ga_3 \subset \Ga_2 \subset \Ga_1 \subset \Ga$ such that 
\begin{enumerate}
\item $\Ga_3$ is a  $p$-group.
\item $\Ga_2/\Ga_3$ is central in $\Ga_1/\Ga_3$\footnote{This claim is not a part of the statement of theorem 0.2 of \cite{LP}. However, $\Ga_1$, $\Ga_2$, and $\Ga_3$ are constructed as the intersection of $\Ga$ with a connected algebraic group, its radical, and its unipotent radical respectively.}, and its order is prime to $p$.
\item $\Ga_1/\Ga_2$ is a product of simple finite groups of Lie type.
\item The index of $\Ga_1$ in $\Ga$ is bounded by a function of $n$ only.
\end{enumerate}
Let $\widetilde{\Ga_i}$ be the subgroups of $S$ such that $\widetilde{\Ga_i}/G_p^1=\Ga_i$.

By our assumptions, $S/V$ has no non-trivial normal $p$-subgroups. Therefore, $\widetilde{\Ga_3}=V$.

If $p$ is large enough, then every $p$-element in $\Ga$ has to be contained in $\Ga_1$. Therefore $\Ga^+\subset \Ga_1$. We clearly have that $\Ga^+/(\Ga_2\cap \Ga^+) \subset \Ga_1/\Ga_2$. Since $\Ga_1/\Ga_2$ is generated by its elements of order $p$, and since every such element can be lifted to an element of $\Ga_1$, we get that, in fact, $\Ga^+/(\Ga_2\cap \Ga^+) = \Ga_1/\Ga_2$. We get a central extension
\[
0 \to \Ga_2/\Ga_3 \to \Ga^+\cdot \Ga_2 /\Ga_3 \to \Ga^+\cdot \Ga_2 /\Ga_2 = \Ga^+/(\Ga^+\cap \Ga_2) = \Ga_1/\Ga_2 \to 0.
\]
This extension splits as a direct product $\Ga^+\cdot \Ga_2/\Ga_3 = P \times A$, where $P$ is a perfect extension and $A$ is abelian. Since the group $P$ is a perfect central extension of a product of finite simple groups of Lie type, and since each finite simple group is generated by its elements of order $p$, we get that $P$ is generated by its elements of order $p$. Therefore we get that $\Ga^+/\Ga_3 \subset P$. Since the group $A$ is contained in $\Ga_2/\Ga_3$, and the order of $\Ga_2/\Ga_3$ is prime to $p$, we get that $\Ga^+/\Ga_3=P$. 

We consider first the extensions of $\rho$ to $S^+\cdot \widetilde{\Ga_2}$. These extensions are governed by a certain element $\be$ in the second cohomology $H^2(S^+\cdot \widetilde{\Ga_2}/V,\bC ^\times)$, as described in Subsection \ref{subsec:ext.of.reps}. By Proposition \ref{prop:extension.p.group} this element has order $p$. By Proposition \ref{prop:cohomology.finite.simple} the sizes of the first and second cohomology groups of finite simple groups are bounded independently of $p$. By Kunneth formula for the cohomology of products, the sizes of $H^1(\Ga_1/\Ga_2,\bC^\times)$ and $H^2(\Ga_1/\Ga_2,\bC^\times)$, and hence of $H^1(P,\bC^\times)$ and $H^2(P,\bC^\times)$, are bounded independently of $p$. Since $A$ is an abelian group and its size is prime to $p$, the sizes of the first and second cohomology groups of $A$ are also prime to $p$. By Kunneth formula again, we get that the size of the second cohomology group $H^2(\Ga^+\cdot \Ga_2/\Ga_3,\bC^\times)$ is prime to $p$, and hence contain no elements of order $p$. Therefore, the representation $\rho$ extends to a representation of $S^+\cdot \widetilde{\Ga_2}$. By Proposition \ref{prop:Isaccs.6.16}, every representation in $\Irr(S^+|\rho)$ is of the form $\chi_0 \otimes \theta$, where $\chi_0$ is a fixed extension, and $\theta$ is a character of $S^+\cdot \widetilde{\Ga_2}/V=P\times A$. We conclude that
\[
\ze_{S^+\cdot \widetilde{\Ga_2}|\rho}(s)=\ze_{P\times A}(s)=\ze_{(S/V)^+}(s)\cdot |A|.
\]
Since $S^+\cdot \widetilde{\Ga_2} =\widetilde{\Ga_1}$ has a bounded index in $S$, we conclude from Lemma \ref{lem:fin.index} that there is a constant $c_1$ such that
\[
c_1^{-1-s}\cdot \ze_{S|\rho}(s) \leq \ze_{S^+\cdot \widetilde{\Ga_2}|\rho}(s) \leq c_1^{1+s}\cdot \ze_{S|\rho}(s).
\]
Finally, since 
\[
\frac{|A|}{[S:S^+]} = \frac{[\Ga^+\cdot \Ga_2/\Ga_3:P]}{[S:S^+]} = \frac{[\Ga_1:\Ga^+]}{[\Ga:\Ga^+]}=[\Ga:\Ga_1]^{-1}
\]
is bounded as a function of $n$, we get the conclusion of the theorem.
\end{proof}

The computation of $\ze_{G_p|\rho}$ is, thus, reduced to the computation of the decomposition tree of $\rho$. In the next subsection, we shall show how to construct the decomposition trees for families---both for varying the representation $\rho$, and for varying the prime $p$ as well.

\subsection{The Family of Decomposition Trees}

Let $\Grass$ be the Grassmanian of subspaces of $\fg_R$, considered as a $\cT_{Hvf}$-definable set (see Subsection \ref{subsec:def.set}). We have the tautological bundle $\widetilde{\Grass} \subset \Grass \times \fg_R$ consisting of pairs $(v,A)$ such that $A$ belongs to the subspace $v$ of $\fg$. We consider $\widetilde{\Grass}$ also as a definable set. The condition that $v\in \Grass$ is closed under Lie brackets, is a definable condition; so is the condition that $v\in \Grass$ is a nilpotent Lie subalgebra of $\fg_R$. We denote the definable subset of $\Grass$ that consists of the nilpotent Lie subalgebras of $\fg$ by $\Grass _U$. For every prime $p$, if $v\in \Grass_U(\bM_p)$, then the set
\[
L(v)=\{A\in \fg_p | \ac(A)\in v\}
\]
is a pro-nilpotent Lie subalgebra of $\fg_p$. If $p$ is large enough, then $\exp(L(v))$ is a pro-$p$ subgroup of $G_p$ and the Orbit Method (Theorem \ref{thm:orbit.method}) holds for it. 

Let $\X$ be the definable set defined in \ref{subsec:reps.first.congruence}. Recall that for every $p$ large enough we have constructed a map
\[
\Phi_p:\X(\bM_p)\to (\fg _p ^1)^\vee.
\]

We wish to extend this map to larger subalgebras of $\fg_p$. If $v\in \Grass_U$, we denote by $v^T$ the set of elements in $\fg_R$ whose transpose is contained in $v$. Given a prime $p$, an element $x=(A,z)\in \X(\bM_p)$, a nilpotent Lie subalgebra $v\in \Grass_U(\bM_p)$, and an element $\theta \in v^T$, we define a linear character $\widetilde{\Phi_p}(x,v,\theta) \in L(v)^\vee$ by
\[
\widetilde{\Phi_p}(x,v,\theta) (B) = \exp \left ( \frac{2\pi i}{z} <A,B> \right ) \cdot \exp \left ( \frac{2\pi i}{p} <\theta,\ac(B)> \right ) 
\]
for every $B\in L(v)$.

Let $\Y\subset \X\times \Grass_U \times \fg$ be the definable set consisting of triples $(x,v,\theta)$ such that $\theta \in v^T$. We denote by $\widetilde{\Xi_p}$ the Orbit Method map from $L(v)^\vee$ to $\Irr(\exp(L(v)))$, and denote by $\widetilde{\Psi_p}$ the composition of $\widetilde{\Phi_p}$ and $\widetilde{\Xi_p}$. We get a diagram

\begin{equation} \lbl{diag:X.Y}
\xymatrix{\Y(\bM_p) \ar[rr]^{\widetilde{\Phi _p}} \ar[d] && \bigsqcup _{v\in \Grass_U(\bM_p)} L(v)^\vee \ar[r]^{\widetilde{\Xi_p}} \ar[d] & \bigsqcup_{v\in \Grass_U(\bM_p)} \Irr(\exp(L(v))) \ar[d] \\ \X(\bM_p) \ar[rr]^{\Phi _p} && (\fg^1_p)^\vee \ar[r]^{\Xi_p} & \Irr(G_p^1)}
\end{equation}
where the leftmost vertical arrow is the projection to the first coordinate, and the other two vertical arrows are the restriction maps. It is easy to see that this diagram commutes.

\begin{defn} Let $G_\cO$ be the definable set
\[
g=(g_{i,j})\in G_\cO \qq \iff \qq g\in G \wedge \val(g_{i,j})\geq 0.
\]
We have that $G_\cO(\bM_p)=G_p$.
\end{defn}

\begin{defn} \lbl{defn:G.action.X.Y} We define actions of $G_\cO$ on $\X$ and $\Y$ in the following way:
\begin{enumerate}
\item If $g\in G_\cO$ and $x=(A,z)\in \X$, we define
\[
g\cdot (A,z) = (g^{-1}Ag,z).
\]
\item Let $g\in G_\cO$ and $y=(x,v,\theta)\in \Y$. Let $w\in \Grass_U$ be the subspace $\Ad(\ac(g))v$. There is a unique $\tau \in \left ( \Ad(\ac(g))v\right )^T$ such that for every $A\in v$,
\[
<A,\theta> = <\ac(g)^{-1}A\ac(g),\tau >.
\]
We define $g\cdot (x,v,\theta)=(g\cdot x,w,\tau)$.
\end{enumerate}
\end{defn}

For every prime $p$, the group $G_p$ acts on each vertex of Diagram (\ref{diag:X.Y}): on the vertices of left column $G_p=G_\cO(\bM_p)$ acts via Definition \ref{defn:G.action.X.Y}, on the vertices of the middle column $G_p$ acts by the coadjoint action, and on the vertices of the right column $G_p$ acts by conjugation.

\begin{lem} All arrows in Diagram (\ref{diag:X.Y}) intertwine the different actions of $G_p$.
\end{lem}

\begin{proof} For $\Phi_p$, $\widetilde{\Phi}_p$, and the leftmost and middle vertical arrows, this is a simple computation. For $\Xi_p$, $\widetilde{\Xi_p}$, and the rightmost vertical arrow, this follows from Theorem \ref{thm:orbit.method}.
\end{proof}

\begin{lem} \lbl{lem:def.normalizer} There is a definable family $\cN \subset \Grass_U \times G$ such that for every $v\in \Grass_U$, the fiber $\cN_v$ is the normalizer of the subgroup $\exp(L(v))$ in $G$.
\end{lem}

\begin{lem} \lbl{lem:def.stab} \begin{enumerate}
\item There is a definable set $S^{char}\subset \Y \times G_V$ such that for every $p$ and every $y\in \Y(\bM _p)$, we have that $S^{char}_y(\bM _p)$ is the stabilizer in $N_{G_p}(L(v))$ of the linear character $\widetilde{\Phi _p}(y)$ of $L(v)$.

\item There is a definable set $S^{rep}\subset \Y \times G_V$ such that for every $p$ and every $y\in \Y(\bM _p)$, we have that $S^{rep}_y(\bM _p)$ is the stabilizer in $N_{G_p}(\exp(L(v)))$ of the representation $\widetilde{\Psi _p}(y)$ of $\exp(L(v))$. 
\end{enumerate}
\end{lem}

\begin{thm} \lbl{thm:family.decomposition.trees} There is a natural number $N$ and a sequence of definable families $\cT _i \subset \X\times (\Grass _U \times U_n)^i$, for $i=1,\ldots,N$, such that if we denote the natural projection from $\cT_{i+1}$ to $\cT_i$ by $\pi _i$, then:
\begin{enumerate}
\item $\cT_1 = \X\times \{ 0\} \times \{ 0\}$.
\item For every prime $p$, every $(x,v_1,\theta _1,\ldots ,v_i,\theta _i)\in \cT_i(\bM_p)$, and every $j\leq i$, we have $\theta _j \in v_j^T$.
\item For every prime $p$ and every $(x,v_1,\theta _1,\ldots ,v_i,\theta _i)\in \cT_i(\bM_p)$, the fiber $\pi_i^{-1}(x,v_1,\theta _1,\ldots ,v_i,\theta _i)$ consists of the tuples $(x,v_1,\theta _1,\ldots ,v_i,\theta _i,v,\theta)$ such that $\exp(L(v))$ is the maximal normal $p$ subgroup of the group $\Stab_{N(v_1)\cap \ldots \cap N(v_i)}\widetilde{\Psi _p}(v_i,\theta _i)$, and $\widetilde{\Phi_p}(v,\theta)$ is an extension of $\widetilde{\Phi_p}(v_i,\theta _i)$.
\item For every prime $p$ and every $(x,v_1,\theta _1,\ldots ,v_N,\theta _N)\in \cT_N(\bM_p)$, the maximal normal pro-$p$ subgroup of the group $\Stab_{N(v_1)\cap \ldots \cap N(v_N)}\widetilde{\Psi _p}(v_N,\theta _N)$ is $\exp(L(v_N))$.
\item There is a definable partition of $\cT_N$ into finitely many pieces such that the semisimple hull of the stabilizer is constant along each piece.
\end{enumerate}
\end{thm}

\begin{proof} We construct $\cT_i$ by induction. The set $\cT_1$ is defined by the first requirement. Suppose we have constructed $\cT_i$. By Lemma \ref{lem:def.stab}, there is a definable family $S\subset \cT_i \times G_V$ such that for every prime $p$ and every $t=(x,v_1,\theta _1,\ldots ,v_i,\theta _i)\in \cT_i(\bM_p)$, the set $S_t(\bM_p) \subset G_p$ is the stabilizer of $\widetilde{\Psi _p}(x,v_i,\theta _i)$ in the intersection $N_{G_p}(\exp(L(v_1)))\cap \ldots N_{G_p}(\exp(L(v_i)))$. By Proposition \ref{prop:filtration.subgroups}, there is a definable sub-family $U\subset S$ such that for every $p$ and $t\in \cT_i(\bM_p)$, the fiber $U_t(\bM_p)$ is the maximal normal pro-$p$ subgroup of $S_t(\bM_p)$. Hence we have a definable family $\cT'_i \subset \cT_i \times \Grass$ such that for every $p$ and $t$ as above, $(\cT_i')_t(\bM_p)$ is the required unipotent radical. It is easy to see that defining $\cT_{i+1}$ as consisting of the tuples $(x,v_1,\theta _1,\ldots ,v_i,\theta _i,v,\theta)$ such that $(x,v_1,\theta _1,\ldots ,v_i,\theta _i,v)\in \cT'_i$ and $\theta \in v^T$, results in a definable family satisfying the second and third requirements.

Finally, for every $i$, either $\dim(v_i) < \dim(v_{i+1})$, or the unipotent radical of the stabilizer is equal to $\exp(L(v_i))$. Therefore the map $\cT_{i+1}\to \cT_i$ is an isomorphism for $i>N:=\dim(G)$.

The last claim follows from Proposition \ref{prop:filtration.subgroups}.
\end{proof}

\begin{lem} There is a definable partition of $\Y$ into sets $\Y_0,\ldots ,\Y_N$ such that if $(x,v,\theta)\in \Y_i$, then
\[
\frac{|Ad(L(v))(\widetilde{\Phi_p}(x,v,\theta))|}{|Ad(G_p^1)\Phi _p(x)|} =p^{2i}
\]
and therefore
\[
\frac{\dim \widetilde{\Psi_p}(x,v,\theta)}{\dim\Psi_p(x)} =p^i
\]
\end{lem}

\begin{proof} The family of stabilizers of $\widetilde{\Phi_p}(x,v,\theta)$ is definable, and therefore its reduction modulo $p$ is also definable. By Theorem \ref{thm:LangWeil}, the dimension of the fibers are definable. This gives the first identity.

The second identity follows from the first.
\end{proof}

\subsection{Proof of Theorem \ref{thm:uniformity}}

We assume first that $p$ is large enough. By Lemma \ref{lem:def.stab} there is a definable family $S^{rep}\subset \Y\times G_V$ such that for every $y=(x,v,\theta)\in\Y(\bM_p)$, the fiber $S^{rep}_y$ is the stabilizer in $N_{G_p}(L(v))$ of the representation $\widetilde{\Psi_p}(x,v,\theta) \in \Irr(\exp(L(v)))$. By the same lemma, there is a definable family $S^{char}\subset \Y\times G_V$ such that for every prime $p$ and $y=(x,v,\theta)\in \Y(\bM_p)$, the set $S^{char}_y(\bM_p)$ is the stabilizer in $N_{G_p}(\exp(L(v)))$ of the character $\widetilde{\Phi_p}(x,v,\theta)\in L(v)^\vee$. 

\begin{defn} Let $\cT_i$ be the decomposition tree constructed in Theorem \ref{thm:family.decomposition.trees}. For every $k$, if $p$ is a prime, $\ell=(x,v_1,\theta_1,\ldots,v_k,\theta_k)\in\cT_k(\bM_p)$, and $0< j\leq k$, we denote 
\[
S^{rep}_0(\ell)=G_p; \qq S^{char}_0(\ell)=G_p; \qq \Psi_0(\ell)=\Psi_p(x),
\]
\[
S^{rep}_j(\ell)=S^{rep}_{(x,v_j,\theta_j)}(\bM_p) ;\qq S^{char}_j(\ell)=S^{char}_{(x,v_j,\theta_j)}(\bM_p) ;\qq \Psi_j(\ell)=\widetilde{\Psi_p}(x,v_j,\theta_j),
\]
and define
\[
W_j(\ell)=\left |\left ( \exp(L(v_j))\cap S^{char}_{j-1}(\ell)\right ) (x,v_1,\theta_1,\ldots,v_j,\theta_j)\right |
\]
and
\[
R_k(\ell)=\frac{[G_p:S^{rep}_{k-1}(\ell)]^{-s}}{\prod_{i\leq k-1} [S^{rep}_i(\ell):S^{char}_i(\ell)\cdot \exp(L(v_i))] \cdot W_i(\ell)}\cdot \left (\frac{\dim\Psi_{k-1}(\ell)}{\dim\Psi_0(\ell)}\right )^{-s}.
\]
\end{defn}

\begin{lem} \lbl{lem:sum.over.leaves} Let $\cT_i$ be the decomposition tree constructed in Theorem \ref{thm:family.decomposition.trees}. For every prime $p$, every $x\in\X(\bM_p)$, and every $k\leq N$, 
\[
\ze_{G_p|\Psi_p(x)}(s)=\sum_{\ell \in (\cT_{k+1})_x(\bM_p)}R_k(\ell)\cdot \ze_{S^{rep}_k(\ell)|\Psi_k(\ell)}(s)
\]
\end{lem}

\begin{proof} Fix $p$ and $x\in\X$. We prove the lemma by induction on $k$. The case $k=0$ is trivial.

Suppose we know the claim for $k-1$. Let $\ell \in(\cT_k)_x(\bM_p)$. The representation $\Psi_{k-1}(\ell)$ is a representation of $\exp(L(v_{k-1}))$, and its stabilizer in $S^{rep}_{k-1}(\ell)$ is just $S^{rep}_k(\ell)$. The maximal normal pro-$p$ subgroup of $S^{rep}_k(\ell)$ is, by the construction, $\exp(L(v_k))$. From Lemma \ref{lem:recursive.rel.zeta} we have that
\[
\sum _{\ell\in(\cT_{k})_x(\bM_p)}R_{k-1}(\ell)\ze_{S^{rep}_{k-1}(\ell)|\Psi_{k-1}(\ell)}(s) = 
\]
\[
\sum _{\ell\in(\cT_{k})_x(\bM_p)} R_{k-1}(\ell) \sum _{\tau\in\Irr(\exp(v_k)|\Psi_{k-1}(\ell))} \frac{[S^{rep}_{k-1}(\ell):S^{rep}_k(\ell)]^{-s}}{|\tau^{S^{rep}_k(\ell)}|} \cdot \left ( \frac{\dim \tau}{\dim \Psi_{k-1}(\ell)} \right ) ^{-s} \cdot \ze_{S^{rep}_k(\ell)|\tau}(s)
\]
For every $\ell \in \cT_k(\bM_p)$, the map $\Psi_k:(\cT_{k+1})_\ell(\bM_p) \to \Irr(\exp(L(v_k))|\Psi_{k-1}(\ell))$ is onto. Hence, by Lemma \ref{lem:change.of.variables},
\[
=\sum _{\ell\in(\cT_{k})_x(\bM_p)}\sum_{\eff \in (\cT_{k+1})_{\ell}(\bM_p)} R(\ell) \cdot \frac{[S^{rep}_{k-1}(\ell):S^{rep}_k(\ell)]^{-s}}{|\Psi_k(\eff)^{S^{rep}_k(\ell)}|\cdot|\Psi_k^{-1}(\Psi_k(\eff))\cap(\cT_{k+1})_\ell(\bM_p)|}\left ( \frac{\dim \Psi_k(\eff)}{\dim \Psi_{k-1}(\ell)} \right ) ^{-s} \cdot \ze_{S^{rep}_k(\ell)|\Psi_k(\eff)}(s)
\]
By definition, $S^{rep}_k(\ell)=S^{rep}_k(\eff)$. By Theorem \ref{thm:orbit.method}, $|\Psi_k^{-1}(\Psi_k(\eff))\cap(\cT_{k+1})_{\ell}(\bM_p)|$ is equal to the size of the orbit of $\eff$ under $\exp(L(v_k))\cap S^{char}_{k-1}(\ell)$, which is just $W_k(\eff)$. Finally, the stabilizer of $\Psi_k(\eff)$ in $S^{rep}_k(\ell)$ is equal to $S^{char}_k(\ell)\cdot\exp(L(v_k))$, and therefore $|\Psi_k(\eff)^{S^{rep}_k(\ell)}|=[S^{rep}_k(\ell):S^{char}_k\cdot\exp(L(v_k))]$. The inductive claim follows.
\end{proof}

By Proposition \ref{prop:red.def.fam} and Lemma \ref{lem:def.root.datum}, there is a definable partition of $\cT_N$ into sets $\cT_N ^1,\ldots,\cT_N^{M}$, and for each part there is a root datum $\cD^i$ such that for each prime $p$ and for each $\ell\in\cT^i_N(\bM_p)$, the group $(S^{rep}_{(x,v_N,\theta_N)}(\ell))^{+}/\exp(L(v_N))$ is isomorphic to $H_{\cD^i}(\bF_p)$. By Lemma \ref{lem:zeta.leaves} we get that there is a constant $c$ such that
\[
c^{-1}\cdot\ze_{S^{rep}_N(\ell)|\Psi_N(\ell)}(s) < \ze_{H_{\cD}(\bF_p)}(s)[S^{rep}_N(\ell):(S^{rep}_N(\ell))^+] < c\cdot\ze_{S^{rep}_N(\ell)|\Psi_N(\ell)}(s).
\]
By Proposition \ref{prop:def.polynom}, there is a $V$-function $\cG$ with domain $\cT_N$ and a constant $c$ such that for every $p$ and every $\ell\in\cT_N(\bM_p)$,
\[
c^{-1}\cG^{\bM_p}(\ell)<[S^{rep}_N(\ell):(S^{rep}_N(\ell))^+]<c\cG^{\bM_p}(\ell).
\]
Together with Proposition \ref{prop:zeta.semisimple.finite} and Proposition \ref{prop:def.polynom}, we get that there is a $V$-function $\cF_{\cD}$ on $\cT_N$ such that for all $p$'s large enough,
\begin{equation} \lbl{eq:uniformity.D}
\ze_{S^{rep}_N(\ell)|\Psi_N(\ell)}(s) \sim \cF_\cD^{\bM_p}(s).
\end{equation}

By Lemma \ref{lem:rel.zeta} and Lemma \ref{lem:change.of.variables} we get that
\[
\ze_{G_p}(s) = \sum _{\rho \in \Irr(G_p^1)}\frac{1}{[G_p:\Stab_{G_p}\rho]}(\dim \rho)^{-s}\cdot\ze_{G_p|\rho}(s) = 
\]
\[
= \int _{x\in\X} \frac{1}{\la(\Psi_p^{-1}(\Psi_p(x)))}\cdot \frac{1}{[G_p:\Stab_{G_p}\Psi_p(x)]}(\dim \Psi_p(x))^{-s}\cdot\ze_{G_p|\Psi_p(x)}(s) dx =
\]
by Lemma \ref{lem:sum.over.leaves} and Theorem \ref{thm:jaikin.precise},
\[
= \int_{x\in\X} p^{\phi_1(x)+s\phi_2(x)}\cdot \frac{1}{[G_p:\Stab_{G_p}\Psi_p(x)]}\cdot \sum _{\ell\in(\cT_{N+1})_{x}(\bM_p)} R(\ell)\cdot \ze_{S^{rep}_N(\ell)|\Psi_N(\ell)}(s) =
\]
\[
\int_{\ell\in\cT_{N+1}(\bM_p)} p^{\phi_1(\pi(\ell))+s\phi_2(\pi(\ell))}\cdot \frac{1}{[G_p:\Stab_{G_p}\Psi_p(x)]} \cdot R(\ell)\cdot \ze_{S^{rep}_N(\ell)|\Psi_N(\ell)}(s).
\]

There is a $V$-function $\cF_{stab}$ with domain $\cT_N$ and a constant $c$ such that for all primes $p$ and all $\ell=(x,v,\theta)\in\cT_{N+1}(\bM_p)$,
\begin{equation} \lbl{eq:uniformity.stab}
c^{-1}\cF_{stab}^{\bM_p}(\ell)<\frac{1}{[G_p:\Stab_{G_p}\Psi_p(x)]} = \frac{[\Stab_{G_p}\Psi_p(x):G_p^1]}{[G_p:G_p^1]}<c\cF_{stab}^{\bM_p}(\ell).
\end{equation}

Similarly, using Proposition \ref{prop:def.polynom}, there is a $V$-function $\cF_R$ with domain $\cT_{N+1}$ such that for all $p$ and $\ell\in\cT_{N+1}(\bM_p)$,
\begin{equation} \lbl{eq:uniformity.R}
c^{-1}\cF_R^{\bM_p}(\ell)<R(\ell)<c\cF_R^{\bM_p}(\ell)
\end{equation}

By (\ref{eq:uniformity.D}), (\ref{eq:uniformity.stab}), and (\ref{eq:uniformity.R}), we get that 
\[
\ze_{G_p}(s)\sim \int_{\ell\in\cT_{N+1}(\bM_p)} p^{\phi_1(\pi(\ell))+s\phi_2(\pi(\ell))}\cdot \cF_{Stab}^{\bM_p}(\ell)\cdot \cF_R^{\bM_p}(\ell)\cdot \cF_{D}^{\bM_p}(\ell).
\]
which proves Theorem \ref{thm:uniformity}.

\begin{lem} \lbl{lem:change.of.variables}
\begin{enumerate}
\item Let $A,B$ be finite sets, let $\phi:A\to B$ be an onto function and let $f:B\to \bC$ be any function. Denote the composition of $f$ and $\phi$ by $g$. Then
\[
\sum_{b\in B} f(b) = \sum _{a\in A}\frac{1}{|\phi^{-1}(\phi(a))|}g(a).
\]
\item Let $(A,\mu)$ be a probability space, let $B$ be a countable set, let $\phi:A\to B$ be a measurable set, and let $f:B\to\bC$ be any bounded function. Denote the composition of $f$ and $\phi$ by $g$. Then
\[
\sum_{b\in B}f(b)=\int_{a\in A}\frac{1}{\mu(\phi^{-1}(\phi(a)))}g(a)da.
\]
\end{enumerate}
\end{lem}

\section{Uniformity of the Local Factors II} \lbl{sec:Uniformity.II}

\subsection{Motivic Integration}

\begin{defn} An $R$-function is a pair $(V,W)$ of definable sets in the language of rings such that $W\subset V\times\bA^n$. If $f=(V,W)$ is a function of type (B) and $p$ is a prime we set
\[
f_p(s)=\sum_{a\in V(\bF _p )}|W_a(\bF _p )|^{-s} .
\]
\end{defn}

\begin{thm} \lbl{thm:mot.int} Let $\cF$ be a bounded $V$-function with domain $X$. Then there are integer constants $A_i ,B_i$, and $R$-functions $f_i$, such that for all but finitely many primes $p$,
\[
\int_{X(\bM_p)}\cF _p (x,s) dx = \sum _i (f_i)_p (s)\cdot \prod _j \frac{p^{A_js+B_j}}{1-p^{A_js+B_j}}.
\]
\end{thm}

\begin{proof} It is enough to prove the theorem for $V$-functions that consist of only one triple $(X,f,\cV)$, where $X\subset \{ (x_1,\dots ,x_n)\in\bA_V^n \qq | \qq \val(x_i)\geq 0\}$ is a definable set, $f:X\to \bA_O$ is a definable function such that $f(x)\geq0$, and $\cV \subset\bA_V^n\times\bA_R^m$ is a definable set. By \ref{thm:elim.quan.valued}, \ref{prop:func.to.ordered} and \ref{prop:red.def.fam} we can assume that there are integral polynomials $P_1 (x),\dots ,P_r (x),Q_1(x),Q_2(x)$, a formula $\varphi(\om _1 ,\dots ,\om _r )$ in the language of fields, a formula $\psi (\ga _1 ,\dots \ga _r)$ in the language of ordered groups, a formula $\xi (x_1,\dots ,x_r, \om _1, \dots ,\om _m)$ in the language of fields, and an integer $e$ such that
\begin{enumerate}
\item $X$ is the set defined by the formula $\varphi(\ac(P_1(x)),\dots,\ac(P_r(x)))\wedge\psi(\val(P_1(x)),\dots,\val(P_r(x)))$.
\item $f(x)=\frac{1}{e}(\val(Q_1(x))-\val(Q_2(x)))$.
\item $\cV$ is the set defined by the formula $\xi(\ac(P_1(x)),\dots,\ac(P_r(x)),\om _1 ,\dots ,\om _m)$.
\end{enumerate}

Let $(Y_\bQ,h_\bQ)$ be a resolution of singularities (see Subsection \ref{subsec:resolution.of.singularities}) for the polynomial $\prod_iP_i(x)\cdot Q_1(x)\cdot Q_2(x)$. Note that $Y_\bQ$ has dimension $n$. We denote the irreducible components of $h^{-1}(D)$ by $E_i$, and denote the closure of $Y_\bQ$ inside $\bP_{\bA_\bZ^n}^m$ by $Y_\bZ$. For any $i$, denote the multiplicity of $(E_i)_{red}$ inside $h^{-1}(D)$ by $N_i$, and denote the multiplicity of $(E_i)_{red}$ inside the divisor $h^*(dx_1 \wedge \dots \wedge dx_n)$ by $\nu _i -1$.

Let $\Sigma$ be the finite set of primes $p$ such that $(Y_\bQ,h_\bQ)$ does not have a good reduction modulo $p$. For every $p\not\in\Sigma$ and for every closed point $\overline{a}$ of $Y_{\bF _p}$ (which we identify with the subscheme of $Y_\bZ$ lying above $\Spec(\bF _p)$) there is a natural number $d$, an open neighborhood $U$, regular functions $u,y_1 ,\dots ,y_n$, and natural numbers $N_1 ,\dots ,N_d$, such that
\begin{enumerate}
\item $y_i$ form a system of parameters for $Y_\bZ$ in $U$.
\item $y_i$ is a local equation for one of the divisor $E_{n_i}$ in $U$ for $i\leq d$.
\item $u$ is invertible in $U$.
\item $(\prod P_i \cdot Q_1 \cdot Q_2)\circ h = uy_1^{N_1} \cdots y_d ^{N_d}$
\item $U$ is irreducible and smooth.
\end{enumerate}
By (1), (4), and (5), there are natural numbers $N_{i,j},M_{k,j}$, $i=1,\dots ,r$ , $j=1,\dots ,d$ and $k=1,2$ and regular functions $u_i , v_k$ that are invertible on $U$ such that
\[
P_i \circ h = u_i \prod y_j ^{N_{i,j}}
\]
and
\[
Q_k \circ h = v_k \prod y_j ^{M_{k,j}}.
\]
By compactness, there are finite number of such neighborhoods that cover $Y_{\bZ_\Sigma}=Y_\bZ\times\Spec\bZ_\Sigma$\footnote{$Y_{\bZ_\Sigma}$ is the part of $Y_\bZ$ that lie over the primes not in $\Sigma$.}. Denote them by $U_1 ,\dots ,U_l$. Let $U_i '=U_i \setminus \cup _{j<i}U_j$. We consider $U_j'$ as definable sets.

Let $p$ be a prime that is not contained in $\Sigma$, and fix $i \in \{ 1,\dots ,l\}$. For every $\overline{a}\in U'_i(\bF _p)$ and for every $z\in Y(\bQ _p)$ such that $\ac(z)=\overline{a}$ we have
\[
P_i \circ h (z)= u_i (z)\prod y_j ^{N_{i,j}}(z)
\]
and similarly
\[
Q_k \circ h (z)= v_k (z)\prod y_j ^{M_{k,j}}(z).
\]
The functions $u_i$ are regular and non vanishing. For almost all primes $p$ we have that $\val\left ( \frac{du_i}{dz_j} \right ) \geq 0$. Therefore the angular component of $u_i(z)$ depends only on the reduction of $z$: $\ac(u_i(z))=\overline{u_i}(\ac(z))$. It follows that $h^{-1}(X)$ can be decomposed into definable sets defined by formulas of the form
\[
\varphi '(\ac(z),\ac(y_i(z))) \wedge \psi (\val(y_i(z)))
\]
where $\varphi '$ is a formula in the language of fields and $\psi$ is a formula in the language of ordered groups. Also, we have that in each piece,
\[
f\circ h (z)= \frac{1}{e}\sum (M_{1,i}-M_{2,i})\val(y_i(z))
\]
and
\[
|h^*(dx_1 \wedge \dots \wedge dx_n)| = \prod |y_i|^{\nu _i -1}|dy_1\wedge \dots \wedge dy_n|
\]
For a set $X$, let $1_X$ be the characteristic function of $X$. Similarly, for a formula $\eta(t)$, let $1_\eta$ be the characteristic function of the set $\{ t | \eta(t) \}$. We have
\[
\int _{\bQ _p ^n}1_X (x) p^{-sf(x)} |\cV (\bF _p)_x|^{-s} = \int _{z\in Y (\bQ _p)} 1_{h^{-1}(X)}(z) p^{-s f\circ h (z)} |\cV (\bF _p)_{h(z)}|^{-s} |h^*(dx_1 \wedge \dots \wedge dx_n )| =
\]
\[
= \sum _{i} \sum _{\overline{a}\in U_i '(\bF _p)}\int _{z\in Y(\bQ _p) \wedge \ac(z)=\overline{a}} 1_{\phi (\overline{a},\ac(y_i(-))}(z) 1_{\psi (\val(y_i (-)))}(z) p^{\sum (-\frac{s}{e}(M_{1,i}-M_{2,i})+\nu _i -1)\val(y_i(z))} |\cV (\bF _p)_{h(z)}|^{-s}.
\]
For every $\overline{a}$, The map $(y_i):\{ z\in Y(\bQ _p) | \ac(z)=\overline{a}\} \to (p\bZ _p)^n$ is a measure preserving bijection. Therefore the above sum equals
\begin{equation} \lbl{eq:B.cone}
\sum _{i} \frac{1}{p^{2n}}\left ( \sum _{(\overline{a},\overline{b})\in V_i (\bF _p)} |W_i (\bF _p)_{(\overline{a},\overline{b})}|^{-s} \right ) \left ( \sum _{\ga \in C} p^{-s(\overline{n} \cdot \ga)+(\overline{m} \cdot \ga)}\right )
\end{equation}
Where $V\subset U_i ' \times\bA_R^n$ is defined by the formula
\[
(x,y)\in V \iff x\in U_i ' \wedge \phi (x,y) ,
\]
$W\subset V \times\bA_R^m$ is defined by the formula
\[
(x,y,z)\in W \iff (x,y)\in U_i ' \wedge \xi (\overline{u}(x)\cdot \prod y_j ^{N_{i,j}},z) ,
\]
$C\subset\bA_O^n$ is defined by the formula
\[
\ga \in C \iff \psi (\ga ) \wedge \bigwedge \ga _i >0 ,
\]
and the functionals $\overline{n} ,\overline{m}$ are defined by
\[
\overline{n}_i=(M_{1,i}-M_{2,i})/e
\]
\[
\overline{m}_i=\nu _i -1 .
\]

For every $i$, the first sum in (\ref{eq:B.cone}) is an $R$-function, which we shall denote by $g_i$. By elimination of quantifiers for the value group (Theorem \ref{thm:elim.quan.ordered}), $C$ can be decomposed into sets defined by conditions of the form
\[
\phi (x) \geq 0 \wedge  \textrm{$\psi (x)$ is divisible by $N$} 
\]
where $\phi (x),\psi (x)$ are affine functionals. After a further decomposition, we can assume each of these sets to be intersection of a cone and a coset of $N\bZ ^n$ for some fixed $N$. It is well known that it is possible to further divide these sets and get that each set is of the form
\[
\{ n_1 v_1 + \dots +n_k v_k \qq | \qq n_1,\dots ,n_k \in \bN \}
\]
for some vectors $v_1 ,\dots ,v_k$ (this fact is used in the desingularization theorem for Toric varieties: see \cite[Section 2.6]{Fu}). On each cone we have to sum a geometric series, so the sum is of the form
\[
\prod _j \frac{p^{A_j s+B_j}}{1-p^{A_j s+B_j}},
\]
and so for every $p\not\in\Sigma$,
\[
\int_{X(\bM_p)}\cF_p(x,s)dx=\sum_i(g_i)_p(s)\prod_j\frac{p^{A_{i,j}s+B_{i,j}}}{1-p^{A_{i,j}s+B_{i,j}}}.
\]

If $p\in\Sigma$, we can resolve the singularities of $\prod_iP_i(x)\cdot Q_1(x)\cdot Q_2(x)$ in $\bQ_p[x]$. By similar arguments to the above, we get that there are integers $n_i,m_i,C_{i,j},D_{i,j}$ such that
\[
\int_{X(\bM_p)}\cF_p(x,s)dx=\sum_i n_i(m_i)^{-s}\prod_j\frac{p^{C_{i,j}s+D_{i,j}}}{1-p^{C_{i,j}s+D_{i,j}}}.
\]

Arguing as in the proof of Theorem \ref{thm:uniformity}, we can find $R$-functions $f_i$ and integers $A_{i,j},B_{i,j}$ such that for all $p$,
\[
\int_{X(\bM_p)}\cF_p(x,s)dx=\sum_i(f_i)_p(s)\prod_j\frac{p^{A_{i,j}s+B_{i,j}}}{1-p^{A_{i,j}s+B_{i,j}}}.
\]
\end{proof}

\begin{rem} Instead of using resolution of singularities, it is possible to prove Theorem \ref{thm:mot.int} using the methods of \cite{CL} or \cite{HK}.
\end{rem}

\subsection{Proof of Theorem \ref{thm:rat.absc}}

\begin{thm} \lbl{thm:uniformity.simple} There is a partition of the set of primes into finitely many Artin sets, and for each Artin set the following is true: There is a finite set $I$, and for each $i\in I$ there are nonnegative integers, $d_i$ and $e_i$, and two finite sequences of non-negative integers, $A_{i,j}$ and $B_{i,j}$, such that the sequence of the functions $\ze _{G_p}(s)-1$ is equivalent to the sequence of the functions
\[
s\mapsto \sum _{i\in I} p^{d_i-e_is} \cdot \prod _j \frac{p^{-A_{i,j}s+B_{i,j}}}{1-p^{-A_{i,j}s+B_{i,j}}}.
\]
Moreover, $e_i +\sum _j A_{i,j} >0$ for every $i$.
\end{thm}

\begin{proof} By Theorem \ref{thm:uniformity} and Theorem \ref{thm:mot.int}, there is a $V$-function $\cF$ and $R$-functions $f_j$ such that
\begin{equation} \lbl{eq:uniformity.simple.1}
\ze _{G_p}(s)-1 \sim \int _{X(\bM _p)}\cF _p (x,s) dx = \sum _j (f_j)_p(s)\cdot \prod _k \frac{p^{A_{j,k}s+B_{j,k}}}{1-p^{A_{j,k}s+B_{j,k}}}.
\end{equation}

Let $f=(V,W)$ be an $R$-function. By Corollary \ref{cor:size.def}, there is a constant $c$, a partition of the primes into Artin sets $\cP_i$, and for each $i$ a finite set $D_i\subset \bN\times\bQ_{>0}$, such that for all $p\in\cP_i$ and $a\in V(\bF_p)$,
\begin{equation} \lbl{eq:LW.estimates.R.func}
\left | |W_a(\bF_p)| - \mu p^d \right | < cp^{d-\frac{1}{2}}
\end{equation}
for some $(d,\mu)\in D_i$. Moreover, if we denote by $L_{d,\mu,p}$ the set of elements $a$ in $V(\bF_p)$ such that (\ref{eq:LW.estimates.R.func}) holds, then
\[
\left | |L_{d,\mu,p}|-\nu p^e \right | < cp^{e-\frac{1}{2}}
\]
for some $e=e(d,\mu)\in\bN$ and $\nu=\nu(e,\mu)\in\bQ_{>0}$.

For every $p\in\cP_i$,
\[
f_p(s)=\sum _{a\in V(\bF _p)}|W_a(\bF _p)|^{-s} = \sum _{(d,\mu)\in D_i} \sum _{a\in L_{(d,\mu,p)}}|W_a(\bF _p)|^{-s} \sim 
\]
\[
\sim \sum_{(d,\mu)\in D_i}\nu(d,\mu) \cdot p^{e(d,\mu)} \cdot \left ( \mu \cdot p^{d} \right ) ^{-s} \sim \sum _{(d,\mu)\in D_i} p^{e(d,\mu)-ds}
\]
which together with (\ref{eq:uniformity.simple.1}) implies the theorem.
\end{proof}

We can finally prove Theorem \ref{thm:rat.absc}:

\begin{proof} {\em (of Theorem \ref{thm:rat.absc})} Recall from Section \ref{sec:Euler.Factorization} that there is a finite index subgroup $\De\subset\Ga$, such that the pro-algebraic completion of $\De$ has finite index in the group
\[
\prod _{p\not\in\Sigma}\uG(\bZ_p) \times \uG(\bC).
\]

By Corollary \ref{cor:finite.index.equal.absc}, it suffices to prove that the abscissa of convergence of the Dirichlet series
\[
\prod _{p\not\in\Sigma}\ze_{G_p}(s) \cdot \ze_{\uG(\bC)}(s)
\]
is rational.

By Theorem \ref{thm:uniformity.simple}, there is a partition of the set of primes into finitely many Artin sets $\cA_1,\ldots,\cA_n$, and for each $\cA_i$ there are constants $d_i,e_i,A_{i,j},B_{i,j}$, such that for $p \in \cA _i$,
\[
\ze _{G_p}(s) -1 \sim \sum _{i\in I} p^{d_i-e_is} \cdot \prod _j \frac{p^{-A_{i,j}s+B_{i,j}}}{1-p^{-A_{i,j}s+B_{i,j}}} .
\]

Since the abscissa of convergence of a product of two Dirichlet series is the maximum of the abscissae of convergence of the two series, it is enough to show that the abscissae of convergence of
\begin{enumerate}
\item $\ze _{\uG (\bC)}(s)$, and
\item $\prod _{p\in \cA _i} \ze _{G_p}(s)$ for $i=1,\ldots ,n$
\end{enumerate}
are rational. As noted in Section \ref{sec:Euler.Factorization}, the abscissa of convergence of $\ze _{\uG (\bC)}(s)$ is rational. By Theorem \ref{thm:jaikin}, the abscissa of convergence of each $\ze _{G_p}(s)$ is rational. Thus, if $\cA_i$ is finite, then the abscissa of convergence of $\prod _{p\in \cA _i} \ze _{G_p}(s)$ is rational. We can assume that $\cA_i$ is infinite, and hence has positive analytic density. By Lemma \ref{lem:equiv.Euler}, it is enough to show that the abscissa of convergence of the product

\begin{equation} \lbl{eq:prod}
\prod _{p\in \cA} \left ( 1+p^{d-es} \cdot \prod _j \frac{p^{-A_{j}s+B_{j}}}{1-p^{-A_{j}s+B_{j}}} \right )
\end{equation}
is rational for all subsets $\cA$ of primes with positive density and nonnegative integers $d,e,A_j,B_j$. We shall show that the abscissa of convergence of this product is 
\[
\max \left \{ \frac{\sum _j B_j +d_i +1}{e_i +\sum _j A_j} , \frac{B_j}{A_j} \right \} .
\]

Every factor in the product has a pole at $s=\frac{B_j}{A_j}$, so the abscissa of convergence is greater than the maximum of those expressions. If 
\[
s > \max \left \{ \frac{B_j}{A_j} \right \}
\]
then there is a constant $D>0$ such that for all $p\in \cA$, $1-p^{-A_js+B_j}>D$. Therefore
\[
p^{\left ( \sum _j -A_j -e_i \right ) s+\sum _j B_j +d_j} < p^{d_i} p^{-se_i} \cdot \prod _j \frac{p^{-A_js+B_j}}{1-p^{-A_js+B_j}}  < \frac{1}{D^m} p^{\left ( \sum _j -A_j -e_i \right ) s+\sum _j B_j +d_j}
\]
and we see that if
\[
\left ( \sum _j -A_j -e_i \right ) s+\sum _j B_j +d_j <-1
\]
then the product (\ref{eq:prod}) is greater than the product $\prod _{p\in \cA} (1+\frac{1}{p})$ which diverges, since $\cA$ has positive analytic density.

Similarly, by comparing the product to the (convergent) product $\prod _{p\in \cA} (1+\frac{1}{p^{1+\ep}})$ we see that if 
\[
s >  \frac{\sum _j B_j +d_i +1}{e_i +\sum _j A_j}
\]
then the product (\ref{eq:prod}) converges.
\end{proof}


\vspace{\bigskipamount}

\begin{footnotesize}
\begin{quote}

Nir Avni\\
Einstein Insitute of Mathematics,\\
The Hebrew University of Jerusalem,\\
Edmund Safra Campus, Givat Ram,
Jerusalem 91904, Israel\\
{\tt avni.nir@gmail.com}\\

\end{quote}
\end{footnotesize}
\end{document}